\documentclass[hidelinks,onefignum,onetabnum]{siamart250211}

\usepackage{lipsum}
\usepackage{amsfonts}
\usepackage{graphicx}
\usepackage{epstopdf}
\usepackage{algorithmic}
\ifpdf
  \DeclareGraphicsExtensions{.eps,.pdf,.png,.jpg}
\else
  \DeclareGraphicsExtensions{.eps}
\fi

\newsiamremark{remark}{Remark}
\newsiamremark{hypothesis}{Hypothesis}
\crefname{hypothesis}{Hypothesis}{Hypotheses}
\newsiamthm{claim}{Claim}
\newsiamremark{fact}{Fact}
\crefname{fact}{Fact}{Facts}

\headers{Operator Inference for Elliptic Eigenvalue Problems}{H, Li, J. Sun, and Z. Zhang}

\title{Operator Inference for Elliptic Eigenvalue Problems\thanks{Submitted to the editors \textcolor{black}{DATE}.
\funding{\textcolor{black}{JS was partially supported by NSF Grant DMS-2109949 and SIMONS Foundation Collaboration Grant 711922.  ZZ was partially supported by the National Natural Science Foundation of China  (Project 92470103), the Hong Kong RGC grant (Project 17304324), an R\&D Funding Scheme from the HKU-SCF FinTech Academy, the Outstanding Young Researcher Award of HKU (2020-21), and Seed Funding for Strategic Interdisciplinary Research Scheme 2021/22 (HKU). }}}}

\author{Haoqian Li\thanks{Department of Mathematics, The University of Hong Kong, Pokfulam Road, Hong Kong SAR, China. 
  \email{lihaoqianlhq@connect.hku.hk}.}
\and Jiguang Sun\thanks{Corresponding author. Department of Mathematical Sciences, Michigan Technological University, Houghton, MI 49931, U.S.A. 
  \email{jiguangs@mtu.edu}.}
\and Zhiwen Zhang\thanks{Corresponding author. Department of Mathematics, The University of Hong Kong, Pokfulam Road, Hong Kong SAR, China. AND Materials Innovation Institute for Life Sciences and Energy (MILES), HKU-SIRI, Shenzhen, P.\,R. China.
  \email{zhangzw@hku.hk}}}

\usepackage{amsopn}

\ifpdf
\hypersetup{
	pdftitle={Operator Inference for Elliptic Eigenvalue Problems},
	pdfauthor={H, Li, J. Sun, and Z. Zhang}
}
\fi

\begin{document}
	
	\maketitle
	
	\begin{abstract}
		Eigenvalue problems for elliptic operators play an important role in science and engineering applications, where efficient and accurate numerical computation is essential. In this work, we propose a novel operator inference approach for elliptic eigenvalue problems based on neural network approximations that directly maps computational domains to their associated eigenvalues and eigenfunctions. Motivated by existing neural network architectures and the mathematical characteristics of eigenvalue problems, we represent computational domains as pixelated images and decompose the task into two subtasks: eigenvalue prediction and eigenfunction prediction. For the eigenvalue prediction, we design a convolutional neural network (CNN), while for the eigenfunction prediction, we employ a Fourier Neural Operator (FNO). Additionally, we introduce a critical preprocessing module that integrates domain scaling, detailed boundary pixelization, and main-axis alignment. This preprocessing step not only simplifies the learning task but also enhances the performance of the neural networks. Finally, we present numerical results to demonstrate the effectiveness of the proposed method.
	\end{abstract}
	
	
	
	\section{Introduction}
	Eigenvalue problems for partial differential equations (PDEs) arise in a wide range of disciplines, including mathematics, physics, and engineering. Classical numerical methods—such as the finite element method, finite difference method, and spectral method—have been extensively studied and successfully applied \cite{feit1982solution, kuttler1984eigenvalues, Babuska1991}. While these methods are effective and reliable, they are often computationally intensive, rendering them less suitable for time-sensitive applications such as real-time simulation and inverse problems \cite{pathak2022fourcastnet, nguyen2023climax}.
	
	Recently, data-driven and deep learning approaches have gained increasing attention for the computation of eigenvalue problems. In \cite{han2020solving}, Han et al. construct two deep neural networks to represent eigenfunctions in high-dimensional settings. Li and Ying \cite{li2022semigroup} propose a semigroup-based method using neural networks to solve high-dimensional eigenvalue problems. In \cite{wang2024computing}, Wang and Xie develop a tensor neural network approach to compute multiple eigenpairs. Ji et al. \cite{ji2024deep} introduce a deep Ritz method for approximating multiple elliptic eigenvalues in high dimensions. Ben-Shaul et al. \cite{ben2023deep} design an unsupervised neural network to compute multiple eigenpairs, and Dai et al. \cite{dai2024subspace} propose a subspace method utilizing basis functions generated by neural networks.

	In contrast with the above studies, which treat one eigenvalue problem (the domain is fixed), we are interested in the development of neural operators. Focusing on elliptic eigenvalue problems, we aim to build deep neural networks that take a general domain as input and output eigenvalues and/or eigenfunctions. Such operator learning is crucial for many query scenarios, e.g., shape optimization, inverse spectral problems, and real-time simulations, for which fast and reliable predictions of eigenvalues and eigenfunctions of different domains are necessary.
	
	The first challenge lies in reformulating the eigenvalue problem as an operator learning task. This reformulation can be approached in various ways—for instance, by representing the domain boundary as a continuous function. 
	In this paper, we adopt a neural-network-oriented strategy. Specifically, we choose an operator representation that aligns with existing neural networks known to perform well on similar problems. Leveraging the strengths of Convolutional Neural Networks (CNNs, \cite{lecun1998gradient}) and Fourier Neural Operators (FNOs, \cite{li2020fourier}), we represent both the computational domains and their associated eigenfunctions as pixelated images. This image-based representation explicitly preserves the geometric features of the domains, which is beneficial for learning and generalization.
	
	The second challenge lies in generating a suitable dataset. Our dataset includes both random polygonal domains and smooth domains constructed using random B\'{e}zier curves. To support learning from pixel-based representations, we introduce a preprocessing module designed to address several key issues. First, we scale each domain to fit within the unit square, ensuring consistent size across samples. To account for the rotational and translational invariance of eigenvalues, we apply a main axis alignment (ma) technique to normalize domain orientation. Since pixelization can introduce representation errors, particularly along rapidly varying boundaries, we apply a detailed pixelization (dp) technique that refines boundary representation without increasing image resolution or computational cost. This preprocessing step standardizes the input data and enables the model to focus on intrinsic geometric features, rather than being distracted by variations in scale, position, or orientation. 
	
	The proposed method primarily falls within the framework of neural operator learning \cite{lu2019deeponet,  bhattacharya2021model, wang2022deepparticle}. Since the introduction of Fourier Neural Operators (FNOs), there has been growing interest in employing them for learning mappings between infinite-dimensional function spaces, with demonstrated success in applications such as weather forecasting and inverse problems \cite{khabou2007shape,  thodi2024fourier, pallikarakis2024application}. FNOs excel in learning complex mappings between function spaces due to their ability to capture both local and global features through Fourier transforms. Our proposed method extends this capability by using CNNs and FNOs to respectively predict eigenvalues and eigenfunctions. Specifically, we employ CNNs to learn the mapping from the domain representations (pixelated images) to eigenvalues, and FNOs to learn the mapping from domains to eigenfunctions. This separation of tasks not only enhances the model's interpretability, but also allows for flexible model expansion \cite{Du2023IP}. A key strength of our method is its generalizability to unseen domain geometries, which is particularly useful 
	for shape optimization and real-time design queries, where rapid evaluation of diverse configurations is critical.
	
	
	We conduct extensive experiments to evaluate the robustness, generalization capability, and predictive accuracy of the proposed method. The results demonstrate that our approach achieves high accuracy in both regression and function approximation tasks for the first few eigenpairs. Notably, the relative error in eigenvalue prediction remains stable across different eigenvalue indices, yielding uniformly low errors for the first 20 eigenvalues. In contrast, the relative error in eigenfunction prediction tends to increase with higher eigenvalue indices, reflecting the growing complexity and oscillatory nature of the corresponding eigenfunctions.
	
	The rest of the paper is organized as follows. In Section 2, we discuss the model setup and data generation. We provide details for the main axis alignment and detailed pixelization. Section 3 presents the structures of the CNN and FNO for eigenvalue and eigenfunction prediction, respectively. The loss functions and training of the networks are also discussed. In Section 4, we present various experiments to demonstrate the effectiveness of networks. We end up with conclusions and future work in Section 5.

	\section{Model Setup and Data Generation}
	In this section, we present the modal problem, the choice of the operator, and the generation of dataset. We consider the representative Dirichlet eigenvalue problem in two dimensions.
	Let $\Omega \subset \mathbb R^2$ be a simply connected bounded Lipschitz domain. The  eigenvalue problem is to find $\lambda \in \mathbb R$ and $u \in H^1(\Omega)$ such that 
	\begin{equation}\label{LaplaceEig}
		-\Delta u = \lambda u \quad \text{in } \Omega \quad \text{and} \quad u = 0 \quad \text{on } \partial \Omega.
	\end{equation}
	
	As stated above, we shall represent $\Omega$ by an image such that the shape information is encoded in pixel values, $1$ for pixels inside the domain and $0$ outside. This representation allows for the incorporation of geometric information. To further simplify the task, we shall treat the eigenvalues and the eigenfunctions separately by introducing two subtasks: $\Omega \to \lambda$ and $\Omega \to u$.
	
	The eigenvalues are invariant to both position and rotation of $\Omega$. In addition, they are scale-dependent. If we multiply the $x$ and $y$ coordinates by a factor $k$, the eigenvalues are scaled by a factor of $\frac{1}{k^2}$. Hence, we first scale a domain such that it fits into the unit square $[0, 1]^2$. Taking these properties into consideration, we introduce a preprocessing module to simplify the task. 
	
	\subsection{Domain Generation}
	To generate a domain, we select $n$ random points in $[0, 1]^2$, enforcing a minimum distance threshold $c$ between them. These points are translated so that the centroid is $(0, 0)$ and then reordered according to their principal arguments, denoted by $p_1, p_2, ..., p_n$. Polygons are obtained by connecting adjacent points with line segments
	\begin{equation}\label{non-smooth}
		L_{i, i+1}(t) = t\cdot p_i + (1-t)\cdot p_{i+1}, \quad t \in [0, 1].
	\end{equation}
	
	To obtain a smooth domain, a cubic B\'{e}zier curve segment is generated to connect $p_i$ and $p_{i+1}$. Then $\partial\Omega$ is obtained by concatenating all the B\'{e}zier curve segments. We briefly describe how to construct a cubic B\'{e}zier curve as follows. The angles between each pair of consecutive points are calculated. A weighted average is then computed
	\begin{equation}\label{Thetai}
		\theta_{i, i+1}^{*} = w \cdot \theta_{i-1, i} + (1-w) \cdot \theta_{i, i+1}, \quad i=1, 2, \ldots, n-1,
	\end{equation}
	where $\theta_i$ is the principle argument of $p_i$ and $w = \arctan(\epsilon)/\pi + 0.5$. Here $\epsilon$ is a parameter that controls the ``smoothness" of the curve with $\epsilon=0$ being the smoothest. Two control points between $p_{i}$ and $p_{i+1}$ are generated 
	\begin{equation}\label{ControlPoints}
		\begin{aligned}
			p_i^* &= p_i + r \cdot (\cos(\theta_{i, i+1}^*), \sin(\theta_{i, i+1}^*)), \\
			p_{i+1}^* &= p_{i+1} - r \cdot (\cos(\theta_{i, i+1}^*), \sin(\theta_{i, i+1}^*)),
		\end{aligned}
	\end{equation}
	where the positive number $r \in [0, 1]$ controls the curvature. For $r = 0$, $p_i^*$ and $p_{i+1}^*$ coincide, respectively, with $p_i$ and $p_{i+1}$, and the curve has larger curvature at the control points. Intermediate values of $r$ produce smoother curves, with maximal smoothness when $r$ = 0.5. When increasing further toward $r = 1$, sharp features start to appear near the crossing of the initial and final curve tangents. The points $p_i^*$ and $p_{i+1}^*$ are used to define a cubic B\'{e}zier curve segment connecting $p_i$ and $p_{i+1}$:
	\begin{equation}\label{Bezier}
		B_{i, i+1}(t) = (1-t)^3 p_i + 3 (1-t)^2 t p_i^* + 3 (1-t) t^2 p_{i+1}^* + t^3 p_{i+1}, \quad t \in [0, 1].
	\end{equation}
	
	The above method for generating smooth domain boundaries is inspired by a discussion on Stack Overflow \footnote{https://stackoverflow.com/questions/50731785/create-random-shape-contour-using-matplotlib} and has also been employed in \cite{viquerata2019supervised} for the generation of random shapes.
	
	\subsection{Dataset Generation}
	The eigenvalue problem \eqref{LaplaceEig} has several properties that bring challenges to the design of effective neural networks. Firstly, it is rotationally invariant. Traditional computer vision techniques typically address this by employing data augmentation, where images are rotated by various angles to increase the diversity of the training data. Another approach is to introduce global pooling layers or more complex architectures to capture rotational invariance. However, these methods do not fully address the difficulty.
	
	We choose a different approach that simplifies the problem: main axis alignment, which aligns the domains along their principal axes. This preprocessing step eliminates the rotational invariance by adjusting the domain's orientation before feeding it into the model. Specifically, given a series of points on $\partial \Omega$, denoted by
	$Q = \begin{pmatrix}
		x_0 & x_1&\ldots&x_n \\
		y_0&y_2& \ldots & y_n \\
	\end{pmatrix}^T$, we compute the covariance matrix of $Q$
	\begin{equation}\label{covariance}
		M = \frac{1}{n-1} \sum_{i=1}^n (Q-\bar{Q})^T (Q-\bar{Q}).
	\end{equation}
	Then, we apply the eigenvalue decomposition of $M$
	\begin{equation}\label{EVD}
		M = D_M \cdot v_M,
	\end{equation}
	where the eigenvectors $v_M$ are sorted in ascending order according to the associated eigenvalues (diagonal of $D_M$). The principal axes are given by these eigenvectors. We then rotate the region coordinates about its centroid using $v_M$
	\begin{equation}\label{rotation}
		Q_{\text{rotated}} = (Q - \bar{Q}) \cdot v_M.
	\end{equation}
	
	Secondly, the problem is flipping invariant about $x$ or $y$ axis. To address this issue, we impose an additional constraint by requiring that the $x$-direction standard deviation in the left half of the region is less than or equal to that of the right half, and that the $y$-direction standard deviation in the lower half is less than or equal to that of the upper half. 
	
	Thirdly, we need to represent a domain as a low-resolution pixelated image. The value of a pixel is $0$ if it is completely outside $\partial \Omega$ and $1$ if it is completely inside. The situation for pixels on the boundary is more complicate. The eigenvalues are sensitive to the pixel value on the boundary if the domain is rather narrow. A simple $0/1$ choice inevitably reduces accuracy. 
    
    This issue is particularly subtle and cannot be resolved through conventional hyperparameter tuning, as illustrated by the following example. Consider four rectangular domains $\Omega_i = (0, w_i) \times (0, h)$ with height $h=1$ and width $w_i = \frac{2+\frac{i}{4}}{32}, i=1, 2, 3, 4$. Using $32 \times 32$ pixels, they have the same image but different eigenvalues. The first Dirichlet eigenvalues are
	\begin{equation}
		\lambda_i = \pi^2 \left(\frac{1}{w_i^2}+\frac{1}{h^2}\right), \quad i=1, 2, 3, 4,
	\end{equation}
	which approximately equal to $2006.21$, $1626.91$, $1346.26$, and $1132.81$.
	
	While increasing the image resolution can help, it significantly increases the computational cost. To this end, we employ edge-sensitive pixelization. Specifically, we first generate a high-resolution image and then use average pooling to down-sample the image to the target resolution. The pixel value near the boundary is no longer binary (0/1) but instead a value between 0 and 1, reflecting the proportion of the region contained within the pixel. Although this technique is simple, it significantly improves accuracy without increasing the training cost or the size of the model.
	
\begin{figure}[!htp]
		\centering
		\includegraphics[width=0.95\textwidth]{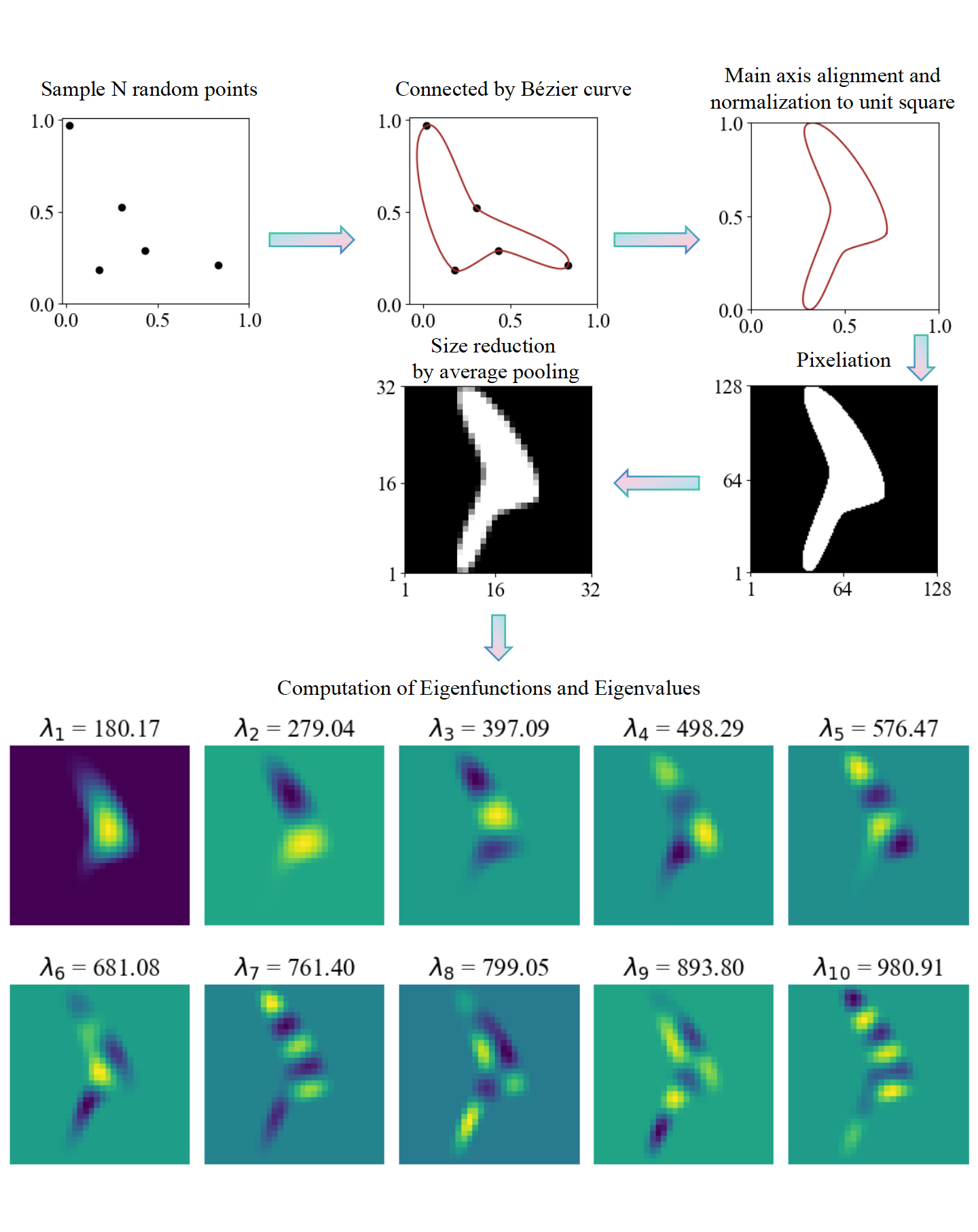}
		\caption{Diagram of dataset generation.}
		\label{diagram_data}
	\end{figure}

\begin{figure}[!htp]
		\centering
		\includegraphics[width=1\textwidth]{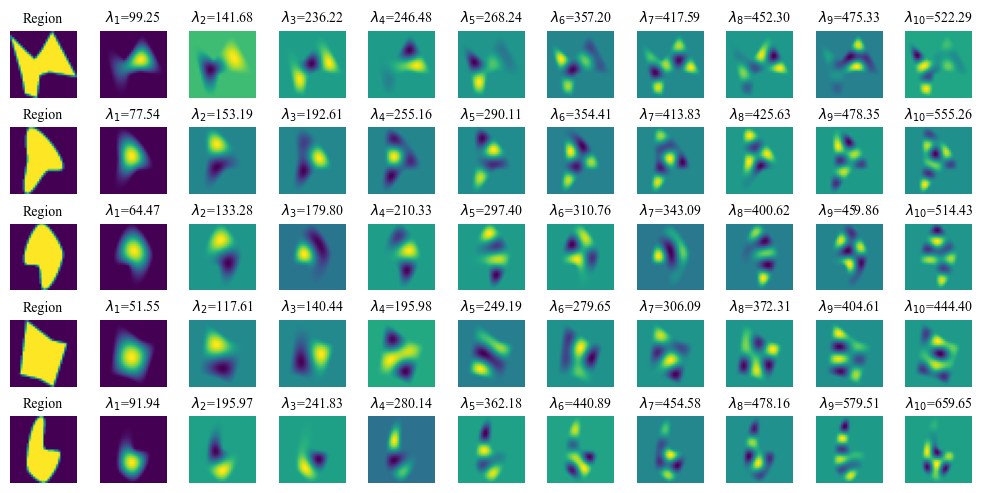}
		\caption{Samples in the dataset. The first 10 eigenvalues and eigenfunctions are shown, as well as their corresponding regions.
		}
		\label{samples_data}
	\end{figure}

	To generate the dateset, i.e., eigenpairs of \eqref{LaplaceEig}, we use a linear Lagrange finite element code on a triangular mesh with mesh size $h \approx 0.01$ (Chp. 3 of \cite{sun2016finite}). The errors of the computed eigenpairs can be ignored for our purposes. Figure \ref{diagram_data} illustrates the workflow of data generation. Figure \ref{samples_data} displays some samples from the dataset. For each domain, it shows the first ten eigenvalues and the associated eigenfunctions.
	
	\section{Methodology}

	\subsection{CNN for Eigenvalue Prediction}
	
	For eigenvalue prediction, we propose a convolutional neural network. CNNs are particularly well-suited for the task due to their powerful ability to capture both global and local geometric features. The local receptive fields of CNNs allow them to effectively extract fine-grained spatial information, while the deeper layers capture global patterns that are critical for understanding the overall shape of the domain.  
	
	In comparison to traditional Deep Neural Networks (DNNs), CNNs exhibit a specialized architecture tailored for processing grid-like data, such as images. One of the primary distinctions lies in how CNNs handle spatial relationships. Unlike DNNs, which treat input data as a flattened vector, CNNs preserve the spatial arrangement of pixels. Furthermore, CNNs reduce the computational burden compared to DNNs by employing weight sharing and pooling operations. The shared weights in convolutional layers drastically reduce the number of parameters compared to fully connected layers in DNNs. This not only enhances the model's ability to generalize but also makes CNNs more computationally efficient. Note that CNNs possess the universal approximation property \cite{he2022approximation}.
	
	The proposed CNN eigenvalue predictor consists of several key components: convolutional layers, max-pooling layers, batch normalization layers, fully connected layers, activation functions (ReLU), and skip connections. Below is a brief overview of each component and its role.
	
	The convolutional layer is the core building block of a CNN. It applies a set of learnable filters (or kernels) to the input data in order to extract local features. For a 2D image input, the convolution operation is given by
	\begin{equation}\label{conv}
		\text{Conv}_{k}(\mathbf{X})(i, j) = (\mathbf{X} * \mathbf{K})(i, j) = \sum_{m=0}^k \sum_{n=0}^k \mathbf{X}(i+m, j+n)\mathbf{K}(m,n),
	\end{equation}
	where $\mathbf{X}\in \mathbb{R}^{d\times d} $ is the input image, $\mathbf{K}$ is the convolutional kernel (filter), $k$ is the size of the filter. 
	
	Max-pooling is applied after a convolutional layer to downsample the spatial dimensions of the feature map. It reduces the computational complexity and helps prevent overfitting. The max-pooling is defined as
	\begin{equation}\label{maxpool}
		\text{Pool}_{k'}(\mathbf{X})(i, j) = \mathop{\max}\limits_{m,n \in [0, k']} \mathbf{X}(i+m, j+n),
	\end{equation}
	where $k'$ is the size of the pooling window. The operation selects the maximum value within each local region, producing a downsampled output.
	
	Batch normalization (Bn) is used to normalize the inputs by adjusting and scaling the activations. This helps to stabilize and speed up the training process. The normalized output is given by
	\begin{equation}\label{batchnorm}
		\text{Bn}(\mathbf{X}) = \gamma \left(\frac{\mathbf{X} - \mu}{\sigma}\right) + \beta,
	\end{equation}
	where $\mu$ and $\sigma$ are the mean and standard deviation of the batch, respectively, and $\gamma$ and $\beta$ are learnable scaling and shifting parameters, respectively.
	
	The fully connected (FC) layer is applied after the convolutional and pooling layers to make predictions based on the extracted features. The FC layer takes the flattened output from the previous layer and computes a weighted sum, followed by a bias,
	\begin{equation}\label{fc}
		\text{FC}_{d_{in}, d_{out}}(\mathbf{x}) = \mathbf{W}_{d_{in}\times d_{out}}\mathbf{x} + \mathbf{b}_{d_{out}},
	\end{equation}
	where $\mathbf{W}_{d_{in}\times d_{out}}$is the weight matrix, $\mathbf{x}\in \mathbb{R}^{d_{in}} $ is the flattened input vector, and $\mathbf{b}_{d_{out}}$ is the bias term.
	
	The activation function introduces non-linearity into the network, enabling it to learn complex patterns. We use the Rectified Linear Unit (ReLU)
	\begin{equation}\label{ReLU}
		\text{ReLU}(x) = \max(0, x).
	\end{equation}
	This function replaces all negative values with zero, introducing sparsity and reducing the risk of vanishing gradients during the training.
	
	Skip connections (or residual connections) allow the model to bypass certain layers and directly pass the output from one layer to a deeper layer. In our case, since we do not need a very deep network, we use skip connections to pass downsampled original input to the deep convolution layers to allow the model to capture global features. These connections also help prevent the vanishing gradient problem and allow for more efficient training. In the context of downsampling, we use a $1 \times 1$ convolution to reduce the spatial dimensions
	\begin{equation}\label{skip}
		\text{Skip}_s(\mathbf{X}) = \mathbf{X} * \mathbf{K'}_{1\times 1}^{(s)},
	\end{equation}
	where $s$ is the stride and $\mathbf{K'}_{1\times 1}^{(s)}$ is a $1\times 1$ conventional kernel with stride $s$, which corresponds to the downsampling factor of the input. This operation effectively reduces the spatial dimensions of the input without changing the depth (number of channels) of the feature maps, making it suitable for skip connections that maintain spatial information while reducing computational cost. The result of this operation is added element-wise to the output of a convolutional layer, combining both local and global information.
	
	For different image sizes, the network needs to be slightly adjusted. For the $32\times 32$ case, we use $3$ convolutional blocks, each containing a convolutional layer with 64/128/256 channels, kernel size $7/5/3$, padding $3/2/1$, and stride $1$. Each convolutional layer is followed by a batchnorm layer, an activation function, and a max-pooling operation with size $2\times 2$, which downsamples the width and height of the feature map by $2$. Next, a skip connection of the downsampled input is added to the output of the max-pooling operation of each convolutional block. The convolutional block can be written as
	\begin{equation}\label{convblock}
		\text{ConvBlock}(\mathbf{X}) = \text{Pool}_{k'}\left(\text{ReLU}(\text{Bn}(\text{Conv}_{k}(\mathbf{X})))\right) + \text{Skip}(\mathbf{X}).
	\end{equation}
	
\begin{figure}[!htbp]
		\centering
		\begin{minipage}[t]{0.49\textwidth}
			\centering
			\includegraphics[width=1\textwidth]{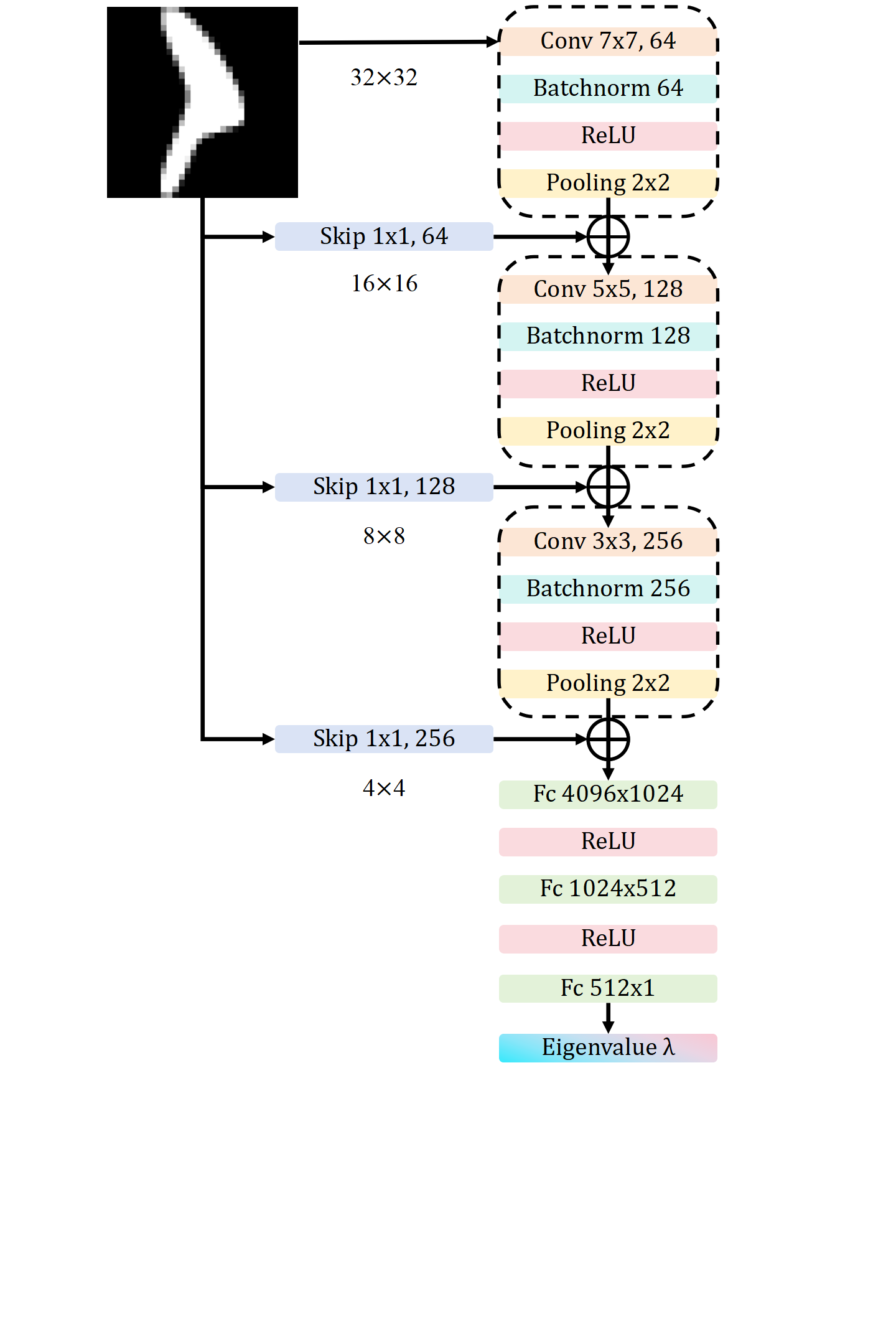}
		\end{minipage}
		\begin{minipage}[t]{0.49\textwidth}
			\centering
			\includegraphics[width=1\textwidth]{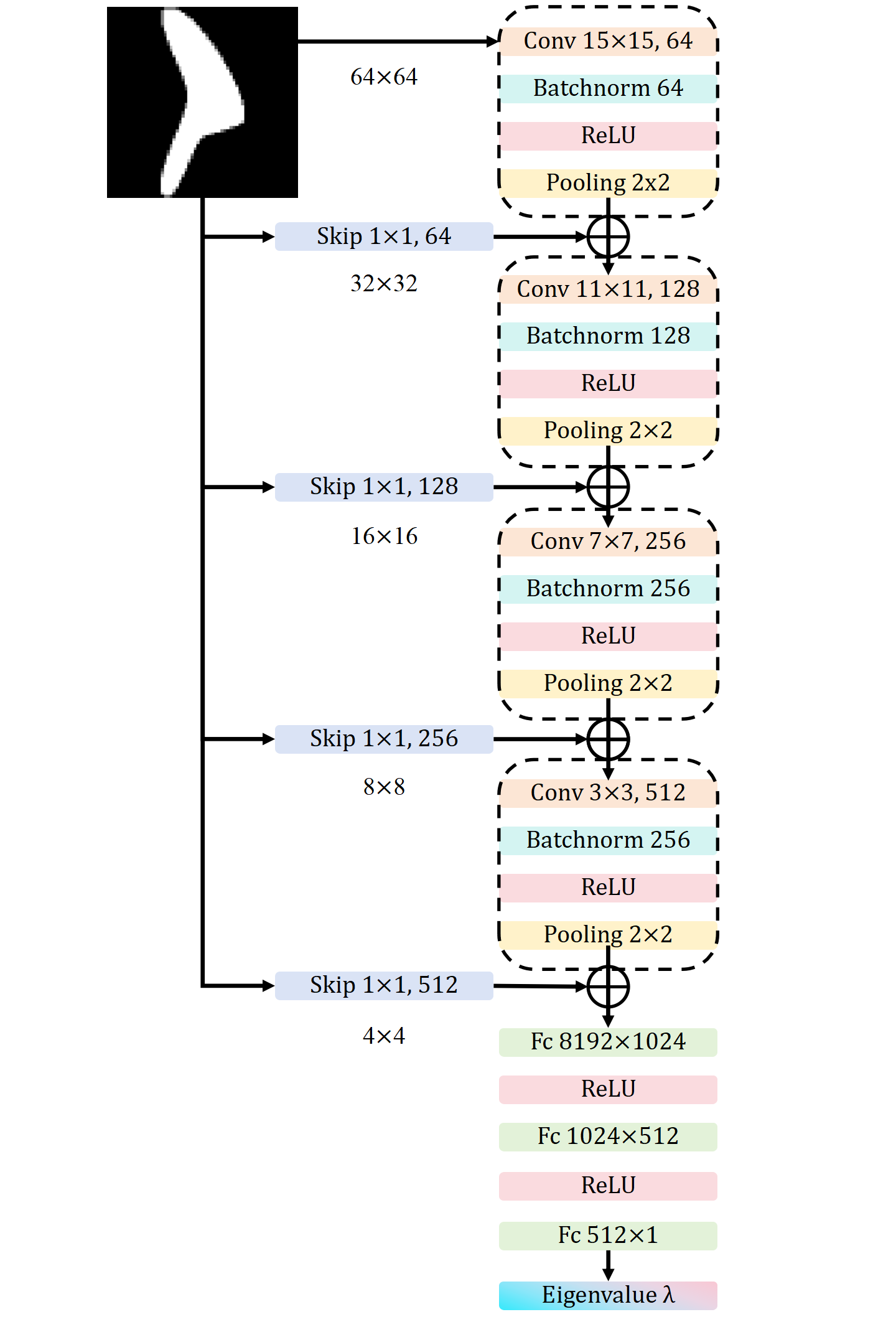}
		\end{minipage}
		\caption{Structures of the proposed CNNs. Left: $32\times 32$ image. Right: $64\times 64$ image.}
		\label{CNN_diagram}
	\end{figure}

	The output of the last convolutional block is then flattened and passed through three fully connected layers and two activation functions in between to make predictions. For $64\times 64$, we simply add another conventional block and adjust the channels, kernel size, and padding accordingly. The rest of the network remains unchanged. The overall structure is illustrated in Figure \ref{CNN_diagram}.
	
	\subsection{FNO for Eigenfunction Prediction}
	We propose a Fourier Neural Operator (FNO) to learn the mapping between images and the eigenfunctions. Domains with complex boundaries or irregular shapes have subtle influences on the eigenfunctions, and FNO's ability to operate in the frequency domain can capture these features efficiently. 
	
	An FNO is organized as a sequence of layers that alternately operate in the spatial domain and the frequency domain. In each Fourier block, the input is first passed through a linear mapping layer $\mathcal{P} $ to lift the dimensionality of the channels. Then, the lifted input passes through multiple Fourier blocks, in which a discrete Fourier transform $\mathcal{F} $ is applied, with a weight tensor multiplication $\mathcal{R} $, which introduces learnable parameters, followed by an inverse Fourier transform $\mathcal{F}^{-1}$. In the Fourier block, the lifted input is also processed in parallel by a local linear transformation $\mathcal{W} $. The output of the above two paths is added element-wise and finally mapped back at the last Fourier block to the output channel dimension via another linear mapping $\mathcal{Q}$. We describe each layer in more detail as follows.
	
	The first layer of the FNO applies a linear mapping $\mathcal{P} $, which operates along the channel dimension of the input. This layer is to elevate the input channels to a higher-dimensional space, enabling the model to learn richer representations of the input data. The transformation is implemented as a fully connected (FC) layer
	\begin{equation}\label{Pmap}
		\mathcal{P}(\mathbf{X}) = \text{FC}_{c_{in}, c_{hd}}(\mathbf{X}),
	\end{equation}
	where $\mathbf{X}\in \mathbb{R}^{c_{in}\times d\times d} $ is the input with $c_in$ channels, $c_{in}$ and $c_{hd}$ are the number of channels in the input and hidden space, respectively. It allows for increased flexibility and capacity in the subsequent layers of the model, as it projects the input features into a higher-dimensional space where more complex interactions can be captured.
	
	Once the input is mapped to a higher-dimensional space, the discrete Fast Fourier Transform (FFT) is applied to compute the frequency components of the input data, which are crucial for understanding periodic or spatially extended features of the domain
	\begin{equation}\label{FFT}
		(\mathcal{F}_M(\mathbf{X}))_c = \sum_{x, y=0}^{d-1} \mathbf{X}_c(x, y) e^{-\frac{2i\pi (M_1x+M_2y)}{d}}, \quad c\in [1, c_{hd}],
	\end{equation}
	where $M=(M_1, M_2)$ are the largest discrete modes for dimensions 1 and 2 (namely $x$ and $y$ direction) and $\mathbf{X}_c$ is the $c$th channel of $\mathbf{X}$. 
    
    FFT transforms the spatial domain data into a frequency domain representation. Then the frequency domain representation is processed by a weight tensor $\mathcal{R}\in \mathbb{C}^{M_{max}\times c_{hd}\times c_{hd}} $ to introduce learnable parameters that enable the model to capture specific patterns or features in the input data
	\begin{equation}\label{Rmap}
		(\mathcal{R} \cdot  \mathcal{F}_M(\mathbf{X}))_{m, c} = \sum_{j=1}^{c_{hd}} \mathcal{R}_{m, c, j} (\mathcal{F}_M(\mathbf{X}))_j, \quad m\in [1, M_{max}], \quad  c\in [1, c_{hd}].
	\end{equation}
	
	Afterward, an Inverse Fast Fourier Transform (IFFT) is applied to map $\tilde{\mathbf{X}}=\mathcal{R} \cdot  \mathcal{F}_M(\mathbf{X})$ to the spatial domain
	\begin{equation}\label{IFFT}
		(\mathcal{F}_M^{-1} (\tilde{\mathbf{X}}))_c = \sum_{\tilde{x}=0}^{d_{\tilde{x}}-1} \sum_{\tilde{y}=0}^{d_{\tilde{y}}-1} \mathbf{\tilde{X}}_c(\tilde{x}, \tilde{y}) e^{2i\pi \left( \frac{(M_1\tilde{x})}{d_{\tilde{x}}} + \frac{(M_2\tilde{y})}{d_{\tilde{y}}} \right) }, \quad c\in [1, c_{hd}].
	\end{equation}
	Finally, the folowing Fourier integral operator $\mathcal{K} $ is applied
	\begin{equation}\label{FourierIntegral}
		\mathcal{K}_M(\mathbf{X})(x, y) = \mathcal{F}^{-1}(\mathcal{R} \cdot  \mathcal{F}_M(\mathbf{X}))(x, y). 
	\end{equation}
	While the data passes through the Fourier integral operator, it also passes through a local linear transformation $\mathcal{W}$ in parallel, which is implemented as a two-stacked $1\times 1$ convolution with an activation function in between. It operates in the spatial domain and allows the model to perform fine-grained adjustments to the features. The transformation $\mathcal{W}$ is defined as
	\begin{equation}\label{Wmap}
		\mathcal{W} (\mathbf{X}) = \text{Conv}_{1\times 1} (\text{GeLU}(\text{Conv}_{1\times 1}(\mathbf{X}))),
	\end{equation}
	where $\text{GeLU}$ is the Gaussian Error Linear Unit
	\begin{equation}\label{GeLU}
		\text{GeLU}(x) = x \cdot \Phi(x), \quad \text{where} \quad \Phi(x) = 0.5x \left(1 +
		\text{tanh}\left(\sqrt{\frac{2}{\pi}} (x + 0.044715 x^3)\right)\right),
	\end{equation}
	a smooth approximation of ReLU that has been shown to improve performance in various deep-learning tasks.

	A Fourier block combining $\mathcal{W}$ and $\mathcal{K}_M(\mathbf{X})$ is defined as
	\begin{equation}\label{FourierBlock}
		\text{FourierBlock}_M(\mathbf{X}) = \text{GeLU}\left( \mathcal{K}_M(\mathbf{X}) + \mathcal{W}(\mathbf{X}) \right).
	\end{equation}
	After several Fourier blocks, the output is passed through a second linear mapping $\mathcal{Q}$, which converts the dimensionality of the channels back to the desired output dimension. This layer is implemented as a fully connected layer and is responsible for mapping the high-dimensional representations learned in the previous layers back to the output space. The mapping $\mathcal{Q}$ is defined as
	\begin{equation}\label{Qmap}
		\mathcal{Q}(\mathbf{X}) = \text{FC}_{c_{hd}, c_{out}}(\mathbf{X}).
	\end{equation}
	
	The complete network representation is as follows:
	\begin{equation}\label{FNO}
		\text{FNO}(\mathbf{X}) = \mathcal{Q} \circ \text{FourierBlock}_M^{(n)} \circ \ldots \circ \text{FourierBlock}_M^{(1)} \circ \mathcal{P} (\mathbf{X}).
	\end{equation}
    
\begin{figure}[!htp]
		\centering
		\includegraphics[width=1.00\textwidth]{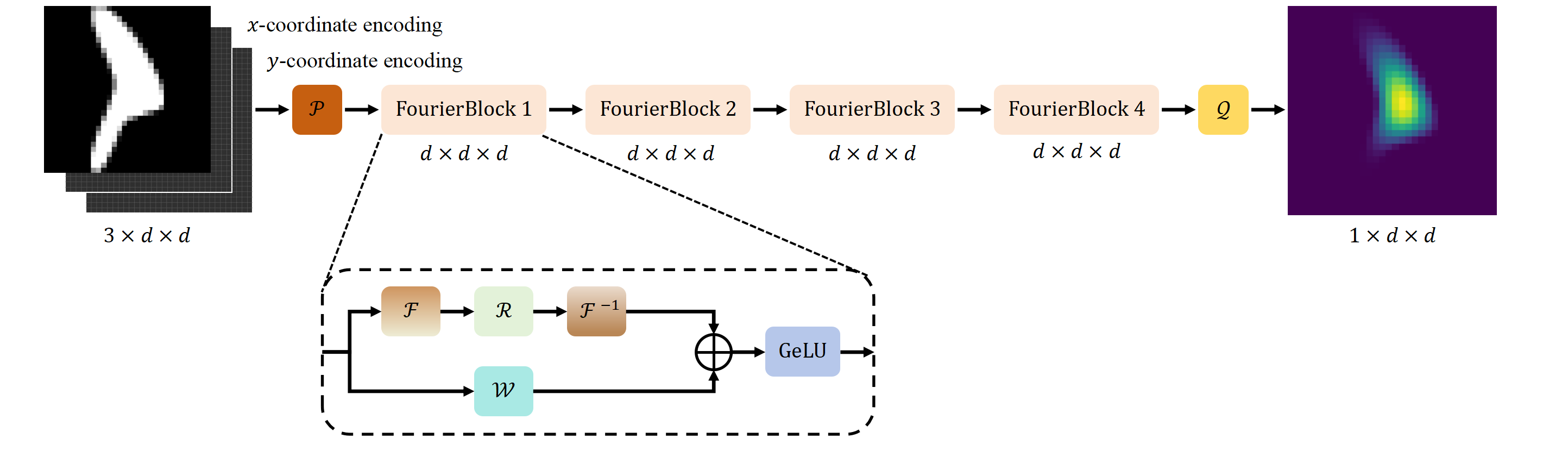}
		\caption{Structure of the proposed FNO.}
		\label{FNO_daigram}
	\end{figure}

	For different input image sizes, the same network architecture is used. Only the number of hidden channels changes. Specifically, for $d\times d$ images, the number of hidden channels is $d$. In the experiments, we use $4$ Fourier blocks and $16$ frequency modes in both directions. To better capture the spatial domain information after the Fourier transform, we incorporate positional encoding into the input. In particular, we augment the original single-channel input by adding $x$ and $y$ coordinate values of each sampling point as additional channels. Consequently, the input image becomes a three-channel representation. The overall architecture of the proposed FNO is illustrated in Figure~\ref{FNO_daigram}.
	
	\subsection{Network Training}
	We describe the loss functions for the networks and the training parameters. For the prediction of the eigenvalue, the loss function is 
	\begin{equation}\label{L1Loss}
		\mathcal{L}_1(\hat{y}, y) = \frac{1}{N} \sum_{i=1}^{N} |\hat{y}_i - y_i|,
	\end{equation}
	where $N$ is the number of data points, $\hat{y}_i$'s are the predicted values, and $y_i$'s are the actual values. The MSE penalizes large deviations between predicted and actual values.
	
	We use the Adam optimizer ($\beta_1 = 0.9, \beta_2=0.999$) for training. The initial learning rate is set to $0.001$, and the learning rate is halved every $10$ steps during training. This allows the model to start with a higher learning rate for faster convergence and gradually reduce it for finer adjustments as the model approaches the optimal solution.
	
	The dataset consists of $20,000$ domain-eigenvalue-eigenfunction triplets, with $10,000$ smooth domains and $10,000$ non-smooth domains. We use $80\%$ of the dataset ($16,000$ input-output pairs) for training and the remaining $20\%$ ($4,000$ input-output pairs) for testing. The batch size is set to $32$. The model is trained for $100$ epochs.
	
	For the prediction of the eigenfunction, the design of the loss function and the training process are more challenging. One needs to take the following key elements into account. (1) Domain of definition. The output image should be $0$ outside $\Omega$. (2) Normalization. It is important that the predicted eigenfunctions are normalized. This prevents the model from learning arbitrary scaling factors. (3) Sign ambiguity. Even if the predicted eigenfunctions are normalized, they may still differ by a sign.
	
	To address the above difficulties, we propose an adaptive relative masked $L_2$ loss
	\begin{equation}\label{ARM2}
		\mathcal{L}_{\text{ARM2}}(\hat{\mathbf{y}}, \mathbf{y}) = \frac{1}{N} \sum_{i=1}^N
		\min \left\{\Vert{\hat{\mathbf{y}}_i \otimes X_i - \mathbf{y}_i}\Vert_2, \Vert{-\hat{\mathbf{y}}_i \otimes X_i - \mathbf{y}_i}\Vert_2 \right\} / \Vert{\mathbf{y}_i}\Vert_2,
	\end{equation}
	where $\otimes$ denotes the Hadamard product (element-wise multiplication), $X_i$ is a $d\times d$ binary mask indicating the domain interior, $\hat{\mathbf{y}}$ and $\mathbf{y}$ are the $d\times d$ predicted and true eigenfunctions, and $\Vert{\,\cdot\,}\Vert_2$ denotes the $L_2$ norm. 
    
    Noted that the output $\mathbf{y}$ is already normalized such that
	\begin{equation}\label{normY}
		\Vert{\mathbf{y}_i}\Vert_2 = C, \quad  i \in \{1, 2, \ldots, N\}.
	\end{equation}

	The loss function \eqref{ARM2} ensures that the value of the predicted eigenfunction within the domain is optimized, while the value outside is ignored. It also takes care of the sign ambiguity by considering both positive and negative versions of the predicted eigenfunction. Finally, the denominator $\Vert{\mathbf{y}_i}\Vert_2$ normalizes the loss by the true eigenfunction norm within the region. It ensures that the predicted eigenfunctions have the same scale, reducing the difficulty of model training.
	
	Again, the Adam optimizer is used. The initial learning rate is $0.001$ and the learning rate decay factor is $0.8$. We set the batch size to be $128$ and train for $50$ epochs. The remaining hyperparameters are the same as the eigenvalue case.
	
	\section{Numerical Experiments}
	We present experimental results, including an ablation study, to show the performance of the proposed networks. Special attention is devoted to the prediction of the first eigenpair, which plays a critical role in applications such as spectral geometry and quantum mechanics. Meanwhile, the successful prediction of other eigenpairs demonstrates the robustness and versatility of our approach.
	
	\subsection{Eigenvalue Prediction}
	To evaluate the performance of the eigenvalue prediction, we use several metrics, including the Root Mean Squared Error (RMSE), R-squared ($R^2$), and Mean Absolute Percentage Error (MAPE).
	
	Let $y_i$ and $\hat{y}_i$ be the true eigenvalues and the predicted eigenvalues, respectively. The RMSE measures the average magnitude of the errors between the predicted and actual values
	\begin{equation}\label{RMSE}
		\text{RMSE} = \sqrt{\frac{1}{N}\sum_{i=1}^N (\hat{y}_i - y_i)^2}.
	\end{equation}
	It is particularly sensitive to large errors, as it squares the differences before averaging. A lower RMSE indicates better prediction accuracy.
	
	Let $\bar{y}$ be the mean of the true eigenvalues. The R-squared metric 
	\begin{equation}\label{R2}
		R^2 = 1 - \frac{\sum_{i=1}^N (\hat{y}_i - y_i)^2}{\sum_{i=1}^N (y_i - \bar{y})^2}
	\end{equation}
	measures the proportion of variance in the true eigenvalues that is explained by the model's predictions. It indicates how well the model captures the underlying relationship between the input and the eigenvalue. An $R^2$ value closer to 1 indicates that the model explains most of the variance in the data, while a value closer to 0 or even smaller than 0 suggests poor performance.
	
	Let $\epsilon$ be a small positive number. The MAPE is defined as
	\begin{equation}\label{MAPE}
		\text{MAPE} = \frac{100\%}{N}\sum_{i=1}^N \left| \frac{\hat{y}_i - y_i}{y_i+\epsilon} \right|,
	\end{equation}
	which measures the average absolute percentage error between the predicted and true eigenvalues. The MAPE is particularly useful for comparing predictions across different scales. 
	
\begin{figure}[!htbp]
		\centering
		\begin{minipage}[t]{0.45\textwidth}
			\centering
			\includegraphics[width=1\textwidth]{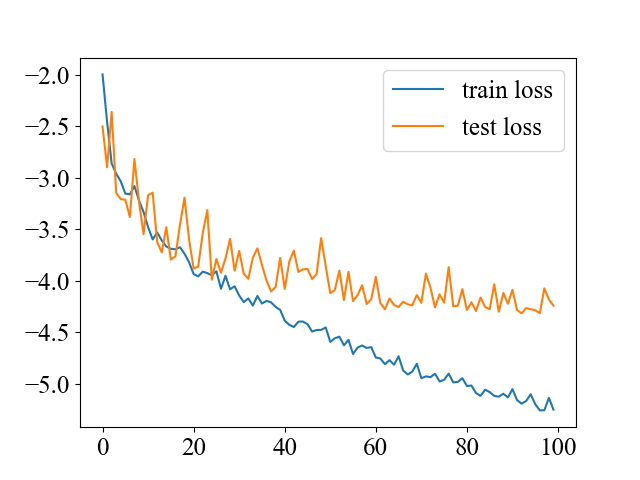}
		\end{minipage}
		\begin{minipage}[t]{0.45\textwidth}
			\centering
			\includegraphics[width=1\textwidth]{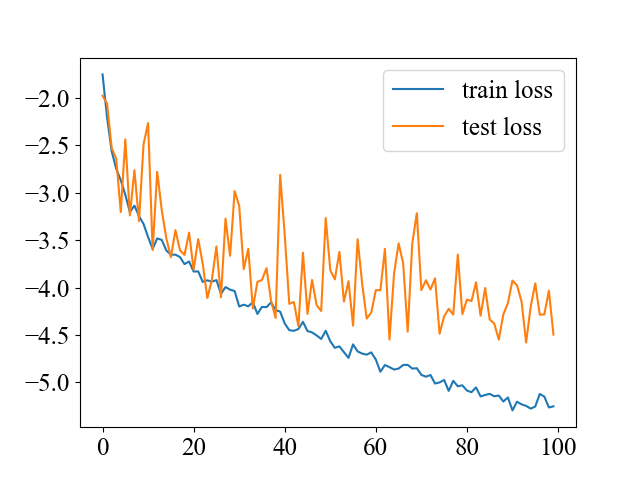}
		\end{minipage}
		\caption{Losses in log scale for the first eigenvalue. Left: $32\times 32$ image. Right: $64\times 64$ image.}
		\label{loss_ev}
	\end{figure}

\begin{figure}[!htbp]
		\centering
		\begin{minipage}[t]{1\textwidth}
			\centering
			\includegraphics[width=1\textwidth]{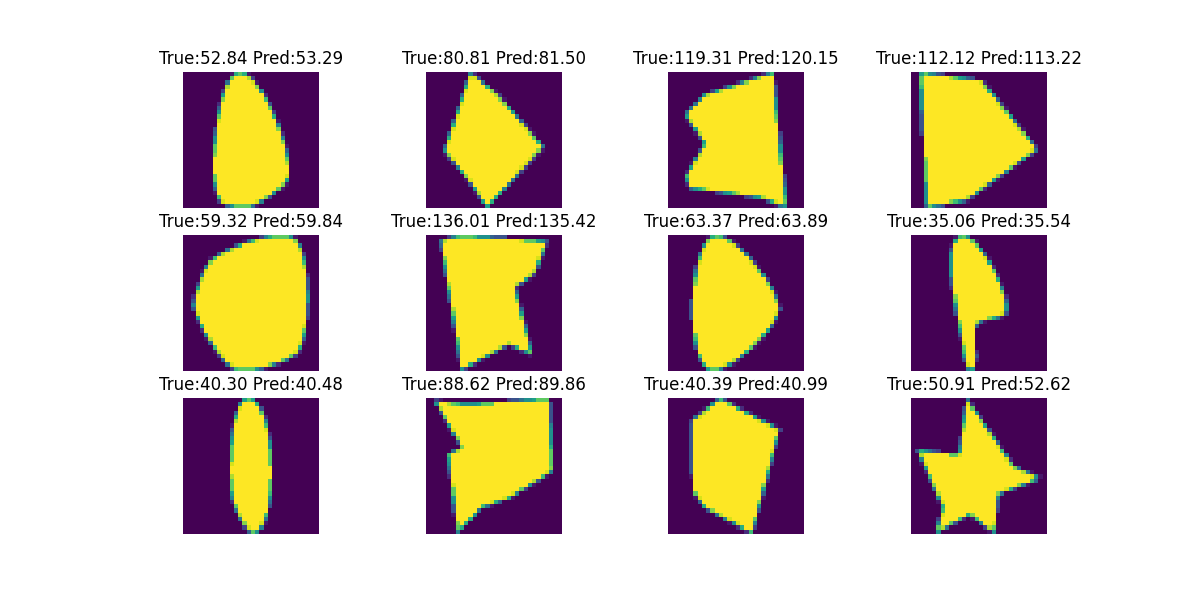}
		\end{minipage}
		\begin{minipage}[t]{1\textwidth}
			\centering
			\includegraphics[width=1\textwidth]{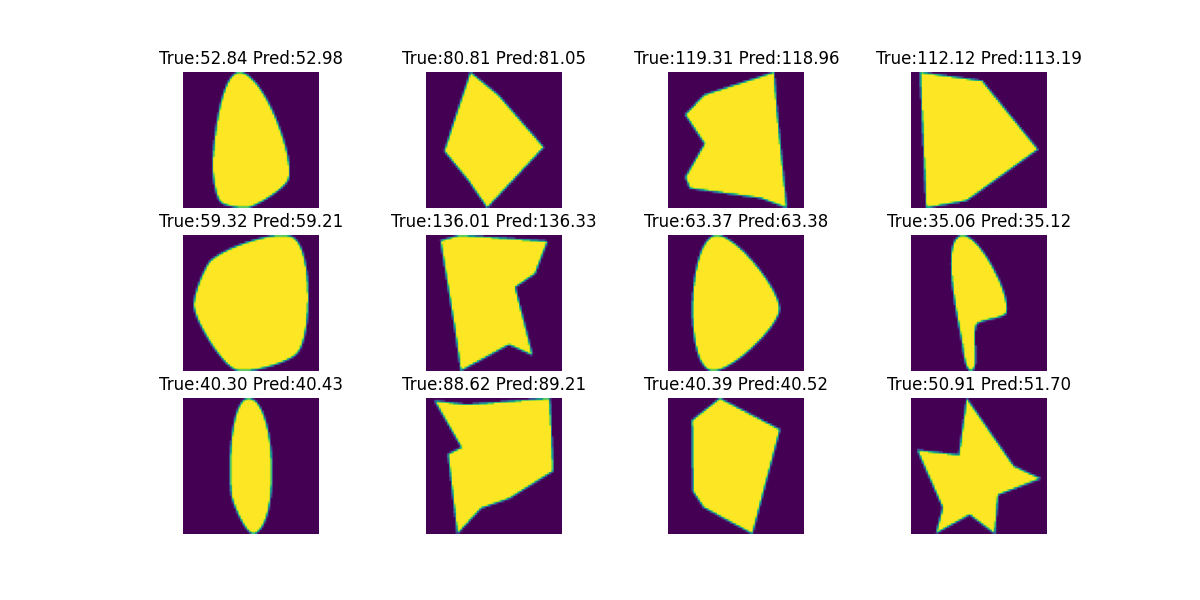}
		\end{minipage}
		\caption{Sample predictions of the first eigenvalue. Top: $32\times 32$. Bottom: $64\times 64$.}
		\label{sample_ev}
	\end{figure}

	Figure~\ref{loss_ev} shows the loss curves for the first eigenvalue prediction for $32 \times 32$ and $64 \times 64$ images in log scale. Both the training and testing losses of the two models stop decreasing significantly around $80$ epochs. Figure~\ref{sample_ev} shows the predictions of the first eigenvalues of some domains, which are close to the true values.
	
	Table~\ref{metric_ev} shows different error metrics for the prediction of the first eigenvalue. The detailed pixelization and main axis alignment significantly improve the model accuracy. The transition from a $32 \times 32$ grid to a $64 \times 64$ yields a slight improvement in the MAPE, from 1.15\% to 0.73\%, and a modest increase in the RMSE. The examples in the rest of the paper use the $32 \times 32$ + dp + ma setting.

	\begin{table}[!htb]
		\footnotesize
		\caption{Performance metrics for the first eigenvalue prediction with different settings. $d\times d$+dp+ma stands for $d\times d$ image with detailed pixelization on the boundaries and main axis alignment. $d\times d$+ma stands for $d\times d$ image with main axis alignment.}
		\begin{center}
			\begin{tabular}{|l|c|c|c|} \hline
				Setting & RMSE & $R^2$ & MAPE \\ \hline
				$32 \times 32$+dp+ma & 2.2121 & 0.9988 & 1.15\% \\
				$32 \times 32$+ma & 5.9933 & 0.9917 & 1.75\%\\
				$32 \times 32$ & 7.3343 & 0.9864 & 1.67\% \\
				$64 \times 64$+dp+ma & 2.2113 & 0.9988 & 0.73\% \\
				$64 \times 64$+ma & 3.0124 & 0.9979 & 1.57\% \\
				$64 \times 64$ & 3.6447 & 0.9966 & 1.78\% \\ \hline
			\end{tabular}
		\end{center}
		\label{metric_ev}
	\end{table}
	
\begin{figure}[!htbp]
		\centering
		\begin{minipage}[t]{1\textwidth}
			\centering
			\includegraphics[width=1\textwidth]{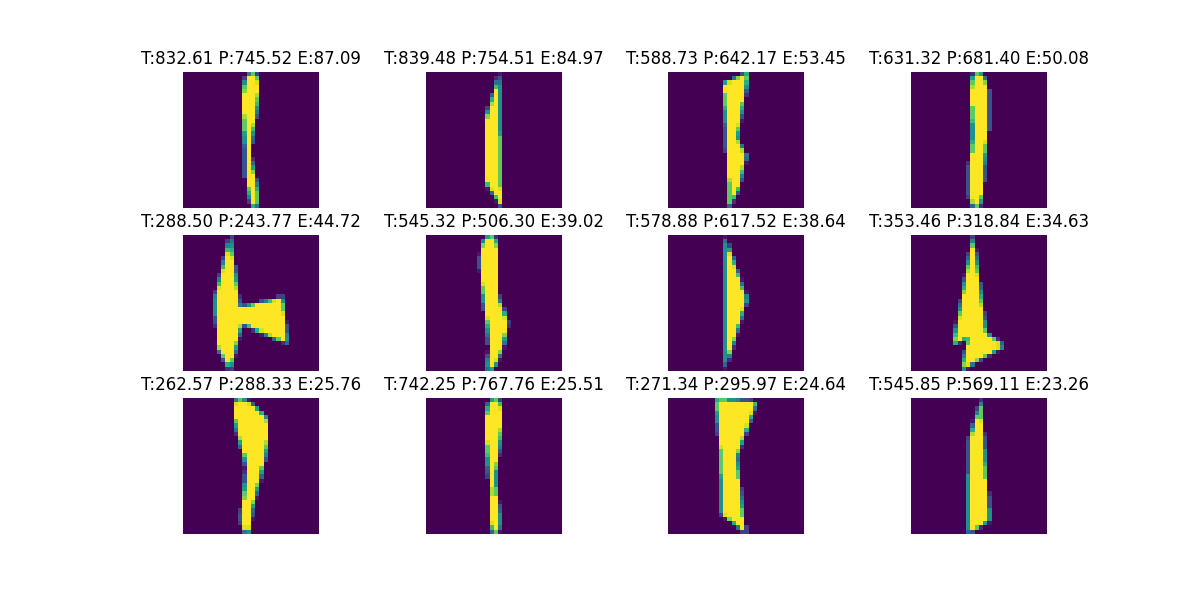}
		\end{minipage}
		\begin{minipage}[t]{1\textwidth}
			\centering
			\includegraphics[width=1\textwidth]{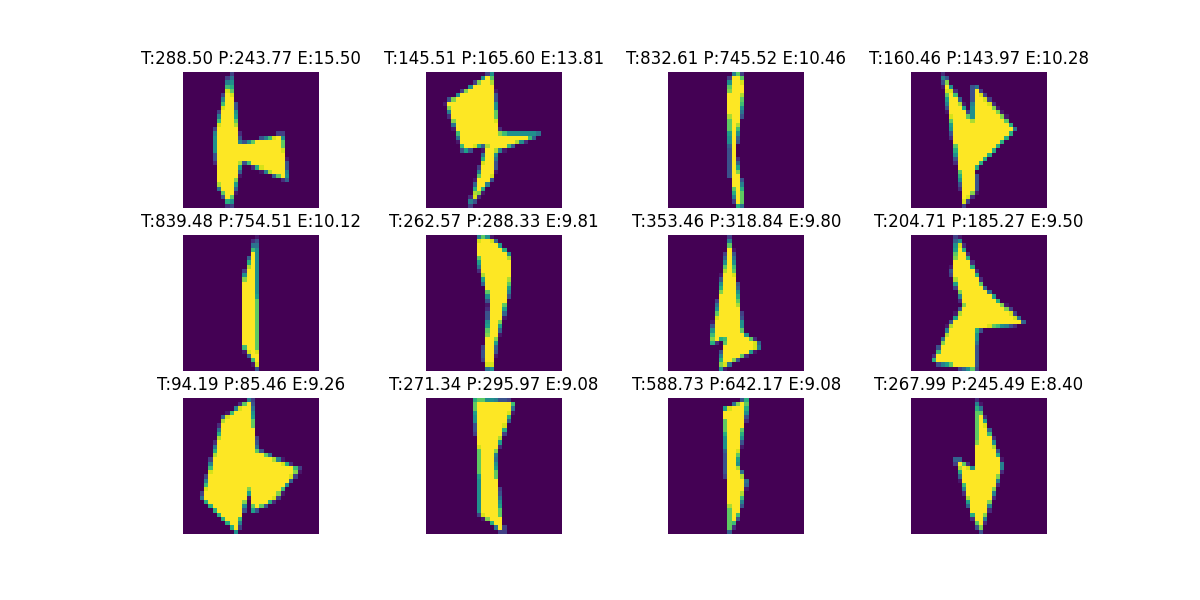}
		\end{minipage}
		\caption{Samples with the largest errors for first eigenvalue predictions. Top: prediction samples with the largest RMSEs. Bottom: prediction samples with the largest MAPEs. "T", "P", and "E" refer to the true, predicted, and error values, respectively. All samples are from the $32 \times 32$+dp+ma setting.}
		\label{large_error_ev}
	\end{figure}

	Figure~\ref{large_error_ev} presents several samples with large errors. They are elongated domains. Such domains are more sensitive to small pixelization errors. More intricate shapes tend to introduce greater approximation errors, leading to higher MAPEs.
	
\begin{figure}[!htbp]
		\centering
		\begin{minipage}[t]{0.32\textwidth}
			\centering
			\includegraphics[width=1\textwidth]{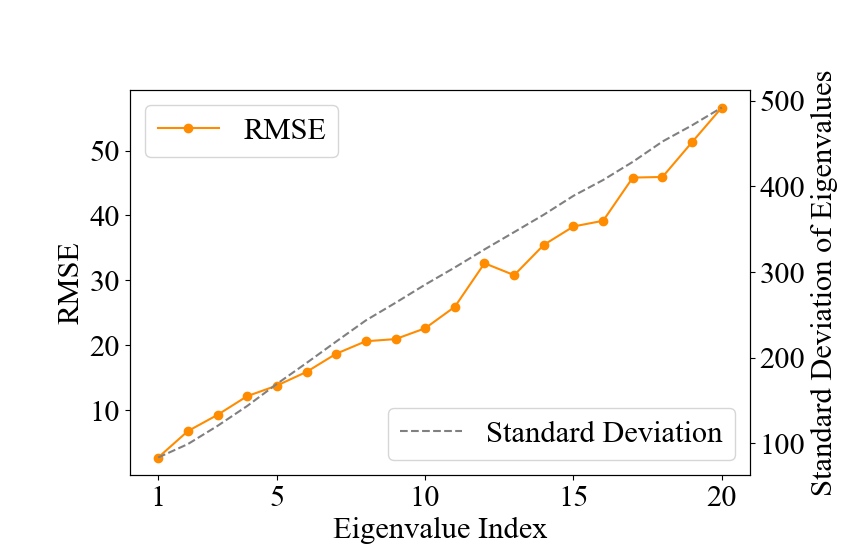}
		\end{minipage}
		\begin{minipage}[t]{0.32\textwidth}
			\centering
			\includegraphics[width=1\textwidth]{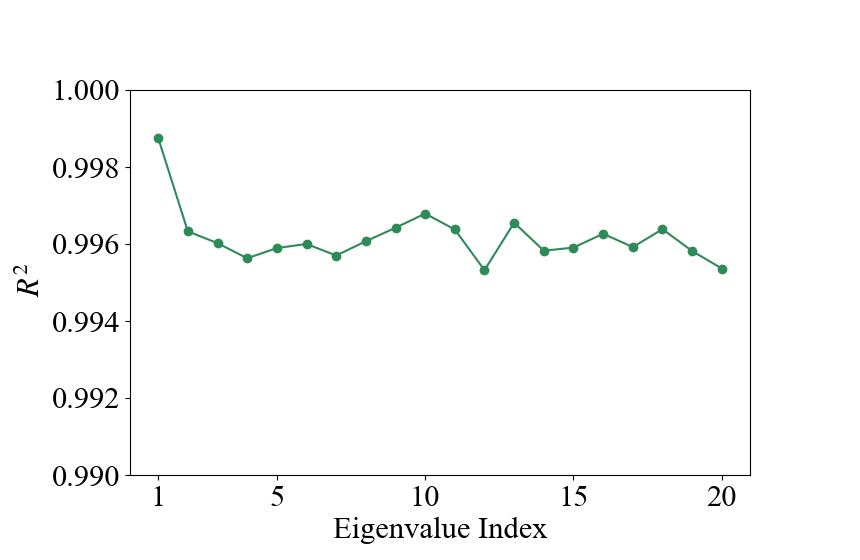}
		\end{minipage}
		\begin{minipage}[t]{0.32\textwidth}
			\centering
			\includegraphics[width=1\textwidth]{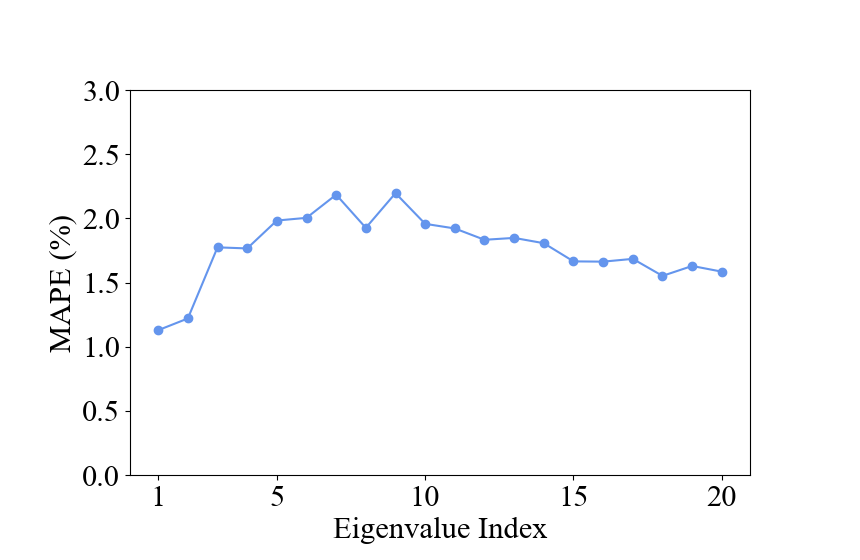}
		\end{minipage}
		\caption{Prediction metrics for the first 20 eigenvalues. Left: RMSE. Middle: $R^2$. Right: MAPE}
		
		\label{metrics_v}
	\end{figure}

	Next, we consider more eigenvalues. Table~\ref{metric_ev_neigs} shows the performance metrics for the 1st, 2nd, and 3rd eigenvalues. It can be seen that RMSE becomes larger for larger eigenvalues, while $R^2$ and MAPE are relatively stable. Figure~\ref{metrics_v} plots RMSE, $R^2$, and MAPE for the first 20 eigenvalues. RMSE increases uniformly with the eigenvalue index, and the growth rate aligns with the standard deviation of the corresponding eigenvalues. $R^2$ and MAPE remain nearly constant.
	
	\begin{table}[!htb]
		\footnotesize
		\caption{Performance metrics for the first three eigenvalues ($32 \times 32$+dp+ma).}
		\begin{center}
			\begin{tabular}{|l|c|c|c|} \hline
				Eigenvalue index & RMSE & $R^2$ & MAPE \\ \hline
				1 & 2.2121 & 0.9988 & 1.15\% \\
				2 & 6.7920 & 0.9963 & 1.22\%\\
				3 & 9.2938 & 0.9850 & 1.78\% \\ \hline
			\end{tabular}
		\end{center}
		\label{metric_ev_neigs}
	\end{table}
	
	\subsection{Eigenfunction Prediction}
	We now present results for the eigenfunction prediction. Three metrics are used: Maximum Absolute Error (MaxAE), Peak Signal-to-Noise Ratio (PSNR), and Relative L1 Error (RelL1). 
	
	The Maximum Absolute Error (MaxAE) is defined as the largest difference between the predicted and true eigenfunctions
	\begin{equation}\label{MaxAE}
		\text{MaxAE} = \max |u_{\text{true}}(x) - u_{\text{pred}}(x)|,
	\end{equation}
	where $u_{\text{pred}}(x)$ and $u_{\text{true}}(x)$ are the predicted and true eigenfunctions evaluated at $x$, respectively. 
	
	The Peak Signal-to-Noise Ratio (PSNR), often used in image and signal processing, is defined as
	\begin{equation}\label{PSNR}
		\text{PSNR} = 20 \log_{10}\left(\frac{\max_x |u_{\text{true}}(x)|}{\rm RMSE}\right),
	\end{equation}
	where $\max_x |u_{\text{true}}(x)|$ is the maximum value of the true eigenfunction and RMSE is the Root Mean Squared Error between the predicted and true functions. 
    
    PSNR provides a measure of the quality of the predicted eigenfunction by comparing the noise (or error) relative to the maximum signal (the true function). A larger PSNR indicates better quality, as the error is smaller relative to the signal's strength. This metric is especially useful for evaluating the overall quality of the predicted eigenfunction, particularly when many small errors may not be captured by MaxAE.
	
	The Relative L1 Error (RelL1) is defined as
	\begin{equation}\label{RelL1}
		\text{RelL1} = \frac{\sum_x |u_{\text{true}}(x) - u_{\text{pred}}(x)|}{\sum_x |u_{\text{true}}(x)|}.
	\end{equation}
	
	The use of MaxAE, PSNR, and RelL1 allows for a comprehensive evaluation of the eigenfunction prediction: MaxAE gives insight into the worst-case errors, PSNR provides a relative measure of the prediction quality, and RelL1 offers an overall picture of the error function. 
	
\begin{figure}[!htbp]
		\centering
		\begin{minipage}[t]{0.45\textwidth}
			\centering
			\includegraphics[width=1\textwidth]{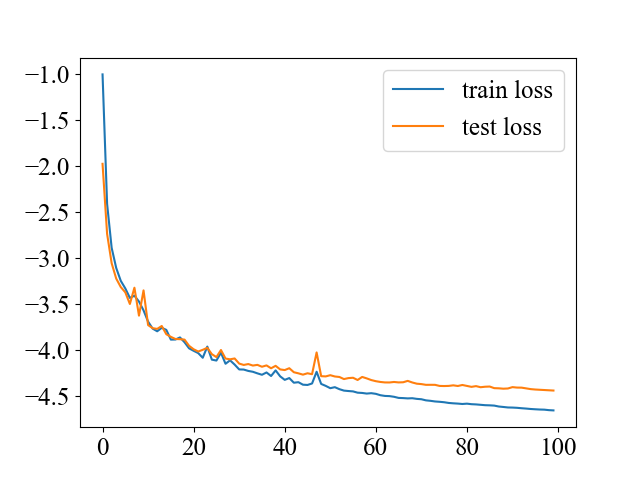}
		\end{minipage}
		\begin{minipage}[t]{0.45\textwidth}
			\centering
			\includegraphics[width=1\textwidth]{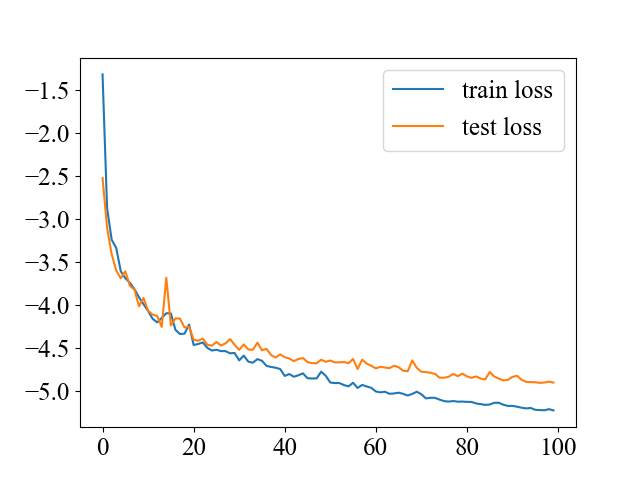}
		\end{minipage}
		\caption{Losses in log scale for the first eigenfunction. Left: $32\times 32$. Right: $64\times 64$.}
		\label{loss_ef}
	\end{figure}

\begin{figure}[!htbp]
		\centering
		\begin{minipage}[t]{1\textwidth}
			\centering
			\includegraphics[width=1\textwidth]{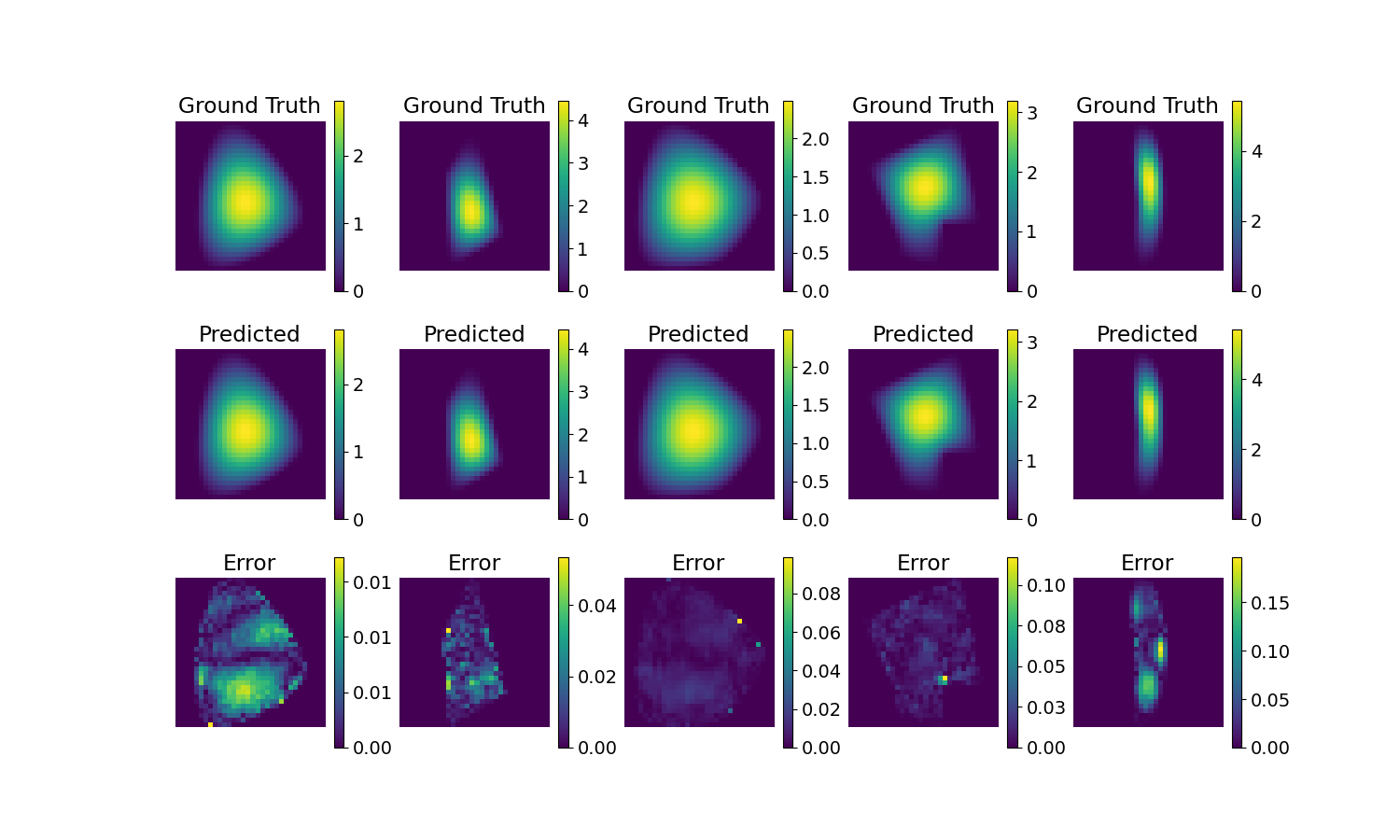}
		\end{minipage}
		\begin{minipage}[t]{1\textwidth}
			\centering
			\includegraphics[width=1\textwidth]{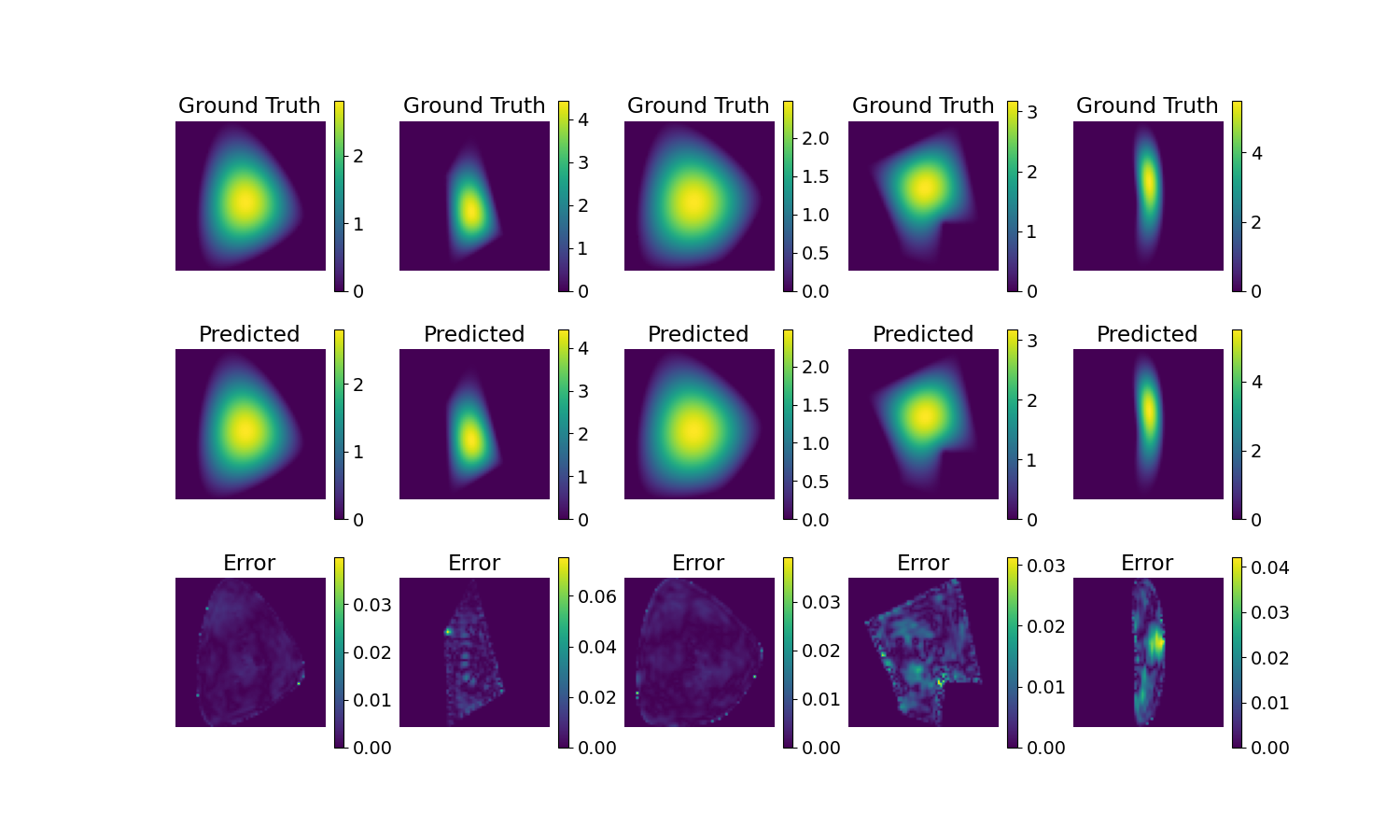}
		\end{minipage}
		\caption{Sample predictions of the first eigenfunction. Top: $32\times 32$. Bottom: $64\times 64$.}
		\label{sample_ef}
	\end{figure}

	Figure~\ref{loss_ef} shows the loss curves for the first eigenfunction prediction on the $32 \times 32$ grid and the $64 \times 64$ grid in log scale. The decreases in training and testing losses slow down significantly after around $40$ epochs, indicating that the models have been sufficiently trained. Figure~\ref{sample_ef} shows sample predictions of the first eigenfunction for both $32 \times 32$ and $64 \times 64$ grids. The predictions are generally close to the actual eigenfunctions, and the scale of errors of both grid sizes is relatively small.
	
	Table \ref{metric_ef} lists MaxAE, PSNR, and RelL1 for the first eigenfunction prediction under different settings. Similar to Table \ref{metric_ev}, the use of detailed pixelization (dp) and main axis alignment (ma) improves the accuracy. For the $32 \times 32$ grid, dp+ma reduces MaxAE from $0.26$ to $0.12$, improves PSNR from $42.64$ to $52.00$, and significantly reduces RelL1 from $3.63\%$ to $1.30\%$. Similar results are observed for the $64 \times 64$ grid. MaxAE decreases from $0.14$ to $0.09$, RelL1 decreases from $1.55\%$ to $0.72\%$, and PSNR increases from $49.88$ to $57.54$.
	
	\begin{table}[!htb]
		\footnotesize
		\caption{Performance metrics for the first eigenfunction prediction with different settings. $d\times d$+dp+ma stands for $d\times d$ grid with detailed pixelization and main axis alignment. $d\times d$+ma stands for $d\times d$ grid with main axis alignment.}
		\begin{center}
			\begin{tabular}{|l|c|c|c|} \hline
				Setting & MaxAE & PSNR & RelL1 \\ \hline
				$32 \times 32$+dp+ma & 0.12 & 52.00 & 1.30\% \\
				$32 \times 32$+ma & 0.27 & 43.14 & 3.57\%\\
				$32 \times 32$ & 0.26 & 42.64 & 3.63\% \\
				$64 \times 64$+dp+ma & 0.09 & 57.54 & 0.72\% \\
				$64 \times 64$+ma & 0.14 & 50.18 & 1.56\% \\
				$64 \times 64$ & 0.14 & 49.88 & 1.55\% \\ \hline
			\end{tabular}
		\end{center}
		\label{metric_ef}
	\end{table}
	
	Table~\ref{metric_ef_neigs} shows the performance metrics for different eigenfunctions. In contrast to Table \ref{metric_ev_neigs}, where the relative error (MAPE) remains stable for different eigenvalues, the prediction of eigenfunctions exhibits a clear trend: both the absolute error (MaxAE) and the relative error (RelL1) increase, and PSNR decreases.
	
	\begin{table}[!htb]
		\footnotesize
		\caption{Performance metrics for the different eigenfunctions ($32 \times 32$+dp+ma).}
		\begin{center}
			\begin{tabular}{|l|c|c|c|} \hline
				Index of Eigenvalue & MaxAE & PSNR & RelL1 \\ \hline
				1 & 0.11 & 53.04 & 1.22\% \\
				2 & 0.37 & 42.13 & 7.04\%\\
				3 & 0.68 & 35.88 & 13.01\% \\ \hline
			\end{tabular}
		\end{center}
		\label{metric_ef_neigs}
	\end{table}
	
\begin{figure}[!htbp]
		\centering
		\begin{minipage}[t]{1\textwidth}
			\centering
			\includegraphics[width=1\textwidth]{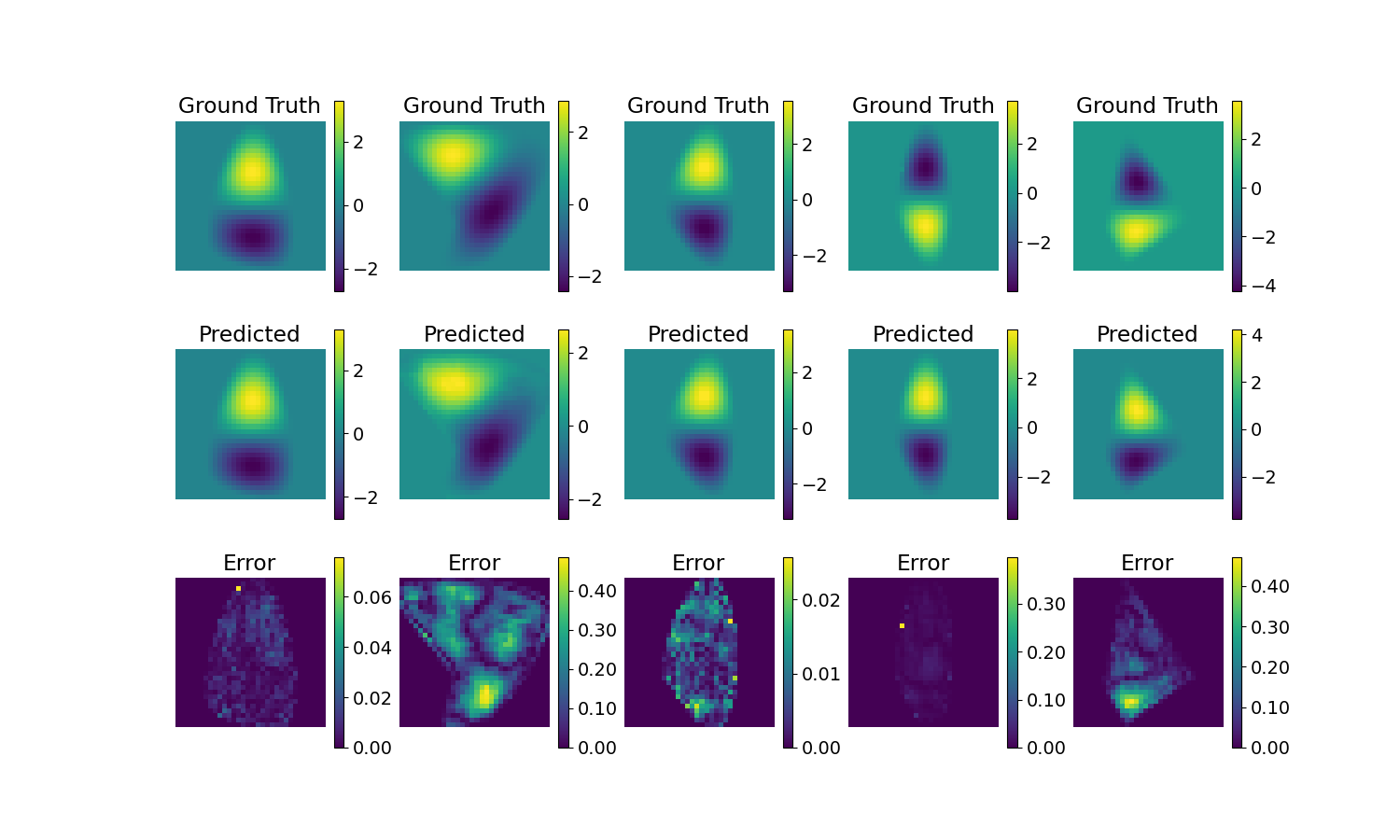}
		\end{minipage}
		\begin{minipage}[t]{1\textwidth}
			\centering
			\includegraphics[width=1\textwidth]{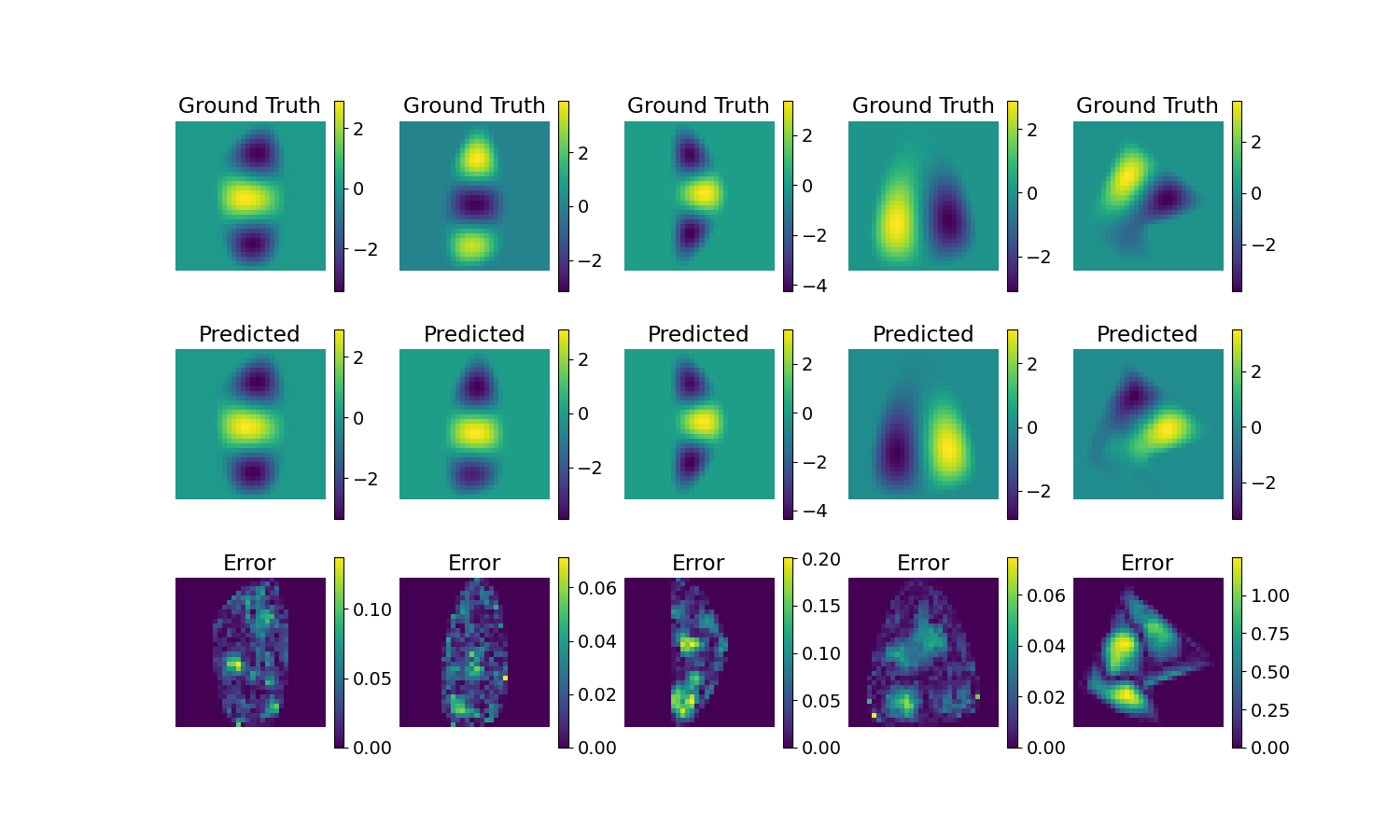}
		\end{minipage}
		\caption{Predictions of the second (top) and third (bottom) eigenfunctions.}
		\label{sample_ef23}
	\end{figure}

	Figure~\ref{sample_ef23} presents several samples of the predicted second and third eigenfunctions. Although the prediction errors are larger than the first eigenfunction, the predicted and true eigenfunctions remain close. For higher eigenfunctions, the model can capture the structure and key features.
	
\begin{figure}[!htbp]
		\centering
		\begin{minipage}[t]{0.32\textwidth}
			\centering
			\includegraphics[width=1\textwidth]{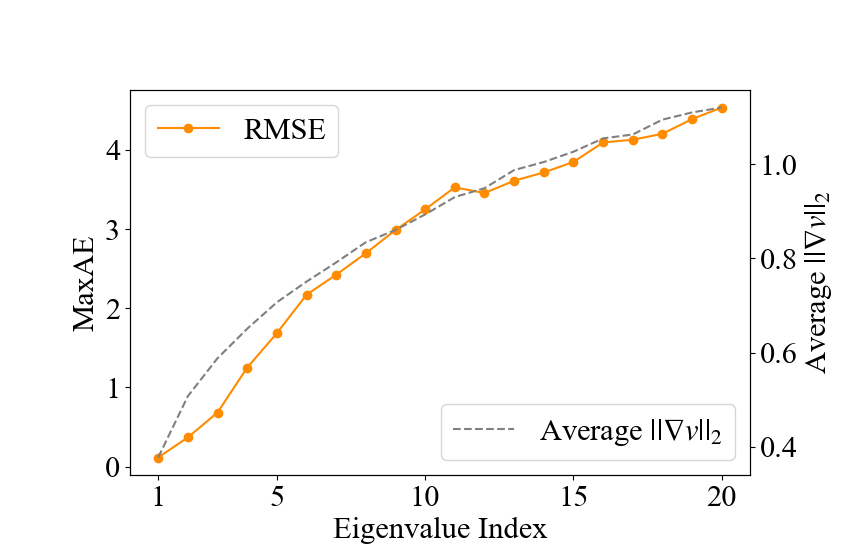}
		\end{minipage}
		\begin{minipage}[t]{0.32\textwidth}
			\centering
			\includegraphics[width=1\textwidth]{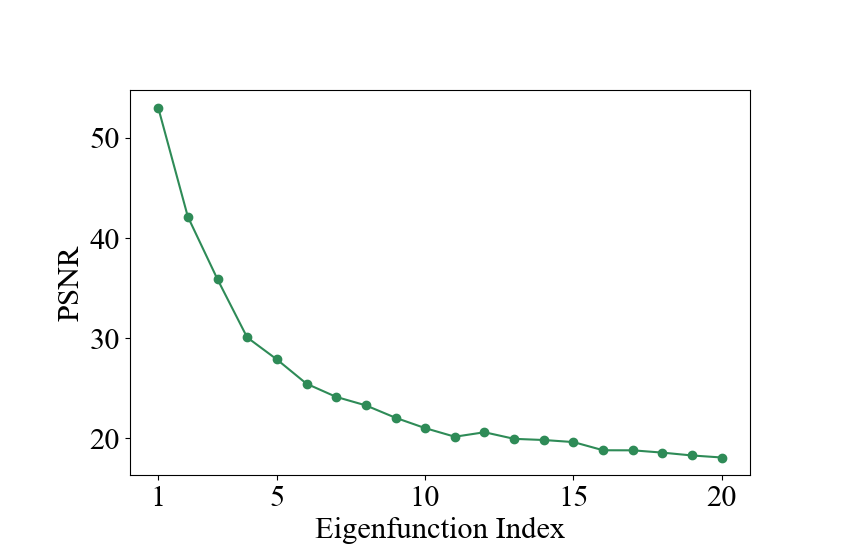}
		\end{minipage}
		\begin{minipage}[t]{0.32\textwidth}
			\centering
			\includegraphics[width=1\textwidth]{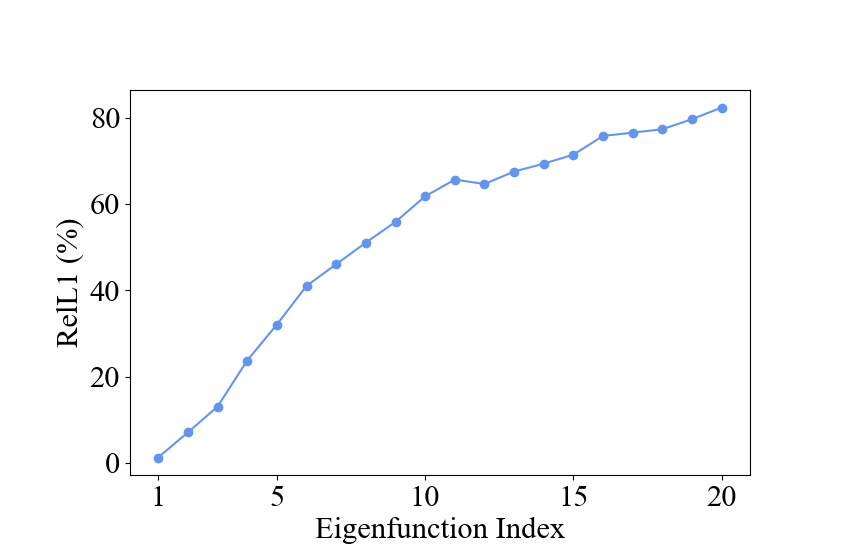}
		\end{minipage}
		\caption{Prediction metrics for the first 20 eigenfunctions. Left: MaxAE. Middle: PSNR. Right: RelL1.}
		\label{metrics_f}
	\end{figure}

	Figure~\ref{metrics_f} shows MaxAE, PSNR, and RelL1 for the first 20 eigenfunctions. The performance declines as the eigenfunction index increases. It is interesting that MaxAE, which reflects the absolute error, aligns with the average gradient norm of the eigenfunctions. This correlation indicates that higher eigenfunctions, with more complex structures and oscillations, pose greater challenges for accurate prediction. The simultaneous decrease in PSNR and increase in RelL1 can partially be attributed to the resolution. 
	
\begin{figure}[!htbp]
		\centering
		\begin{minipage}[t]{1\textwidth}
			\centering
			\includegraphics[width=1\textwidth]{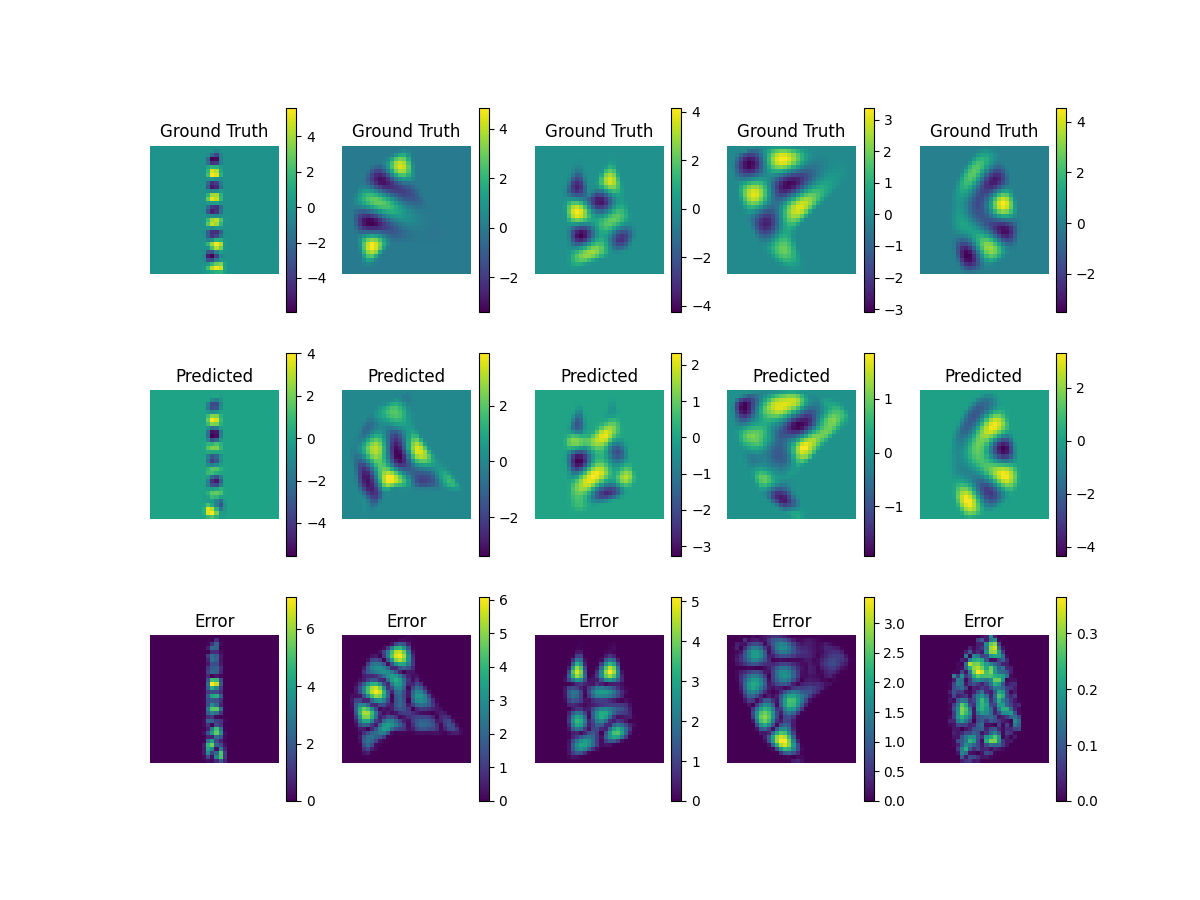}
		\end{minipage}
		\begin{minipage}[t]{1\textwidth}
			\centering
			\includegraphics[width=1\textwidth]{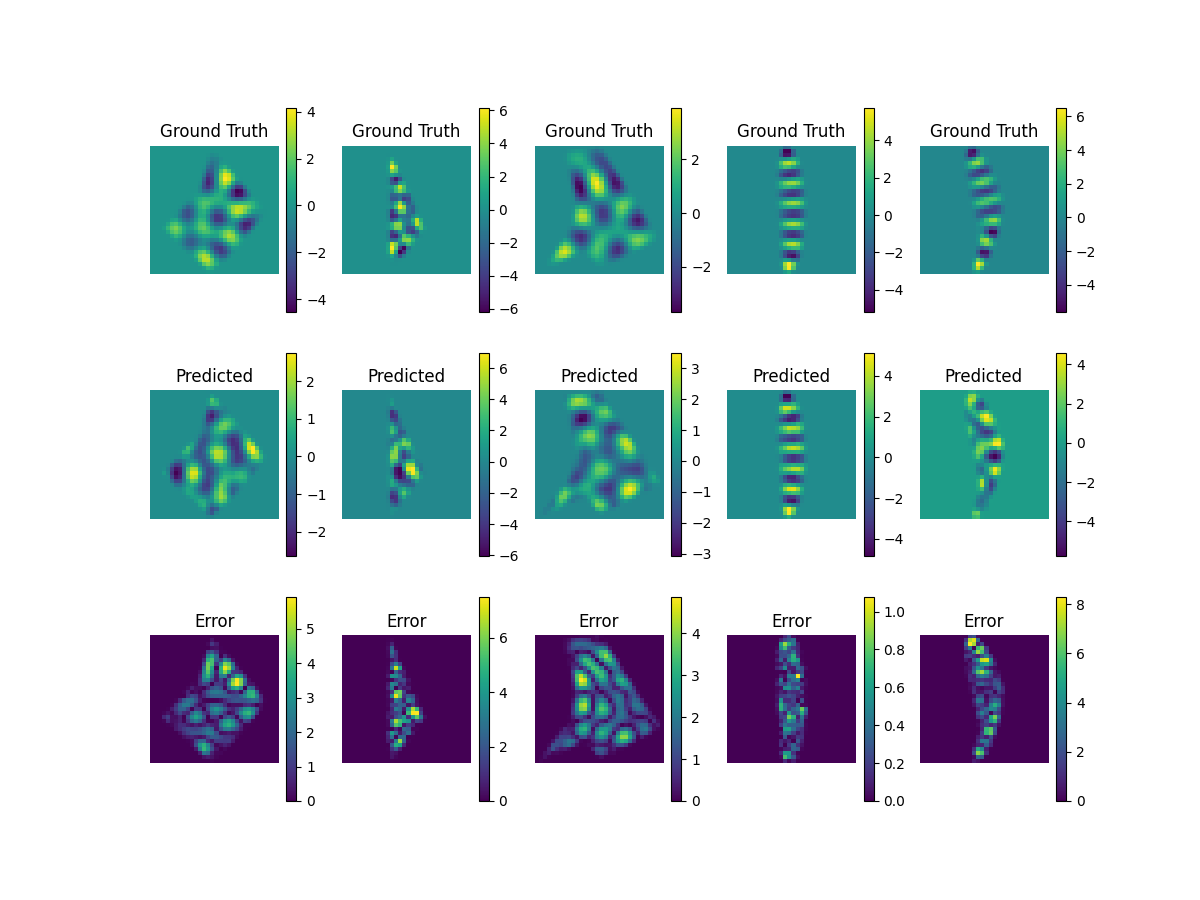}
		\end{minipage}
		\caption{Sample predictions of the 10th and 20th eigenfunctions. Top: the 10th eigenfunction. Bottom: the 20th eigenfunction.}
		\label{sample_ef1020}
	\end{figure}

	Figure~\ref{sample_ef1020} shows predicted results for higher eigenfunctions (10th and 20th). For some samples, the predicted eigenfunctions preserve the overall structures, despite non-negligible errors. Other predicted eigenfunctions differ significantly from the exact ones.
	
	\section{Conclusions}

	In this paper, we proposed a novel operator learning framework for the efficient prediction of eigenvalues and eigenfunctions of elliptic equations on various two-dimensional domains. Unlike most existing neural network methods that focus on single eigenvalue problems or specific high-dimensional settings, our method takes different domains as inputs and outputs the corresponding eigenvalues and eigenfunctions. The key innovations of our approach are threefold. Firstly, we adopt a divide-and-conquer strategy to decompose the task into two subtasks: eigenvalue prediction and eigenfunction prediction.  Secondly, we reformulate the problem in a neural-network-oriented manner, enabling the effective application of CNNs for eigenvalue prediction and FNOs for eigenfunction prediction. Thirdly, we simplify the input problem before passing it into the neural networks, which helps to improve the training efficiency and prediction accuracy. 
	
	To validate the effectiveness of our proposed method, we conducted extensive numerical experiments on a variety of two-dimensional elliptic domains. The results demonstrate that our method achieves significant improvements in terms of prediction accuracy and computational efficiency. In the future, we will extend our approach to more challenging problems, including non-self-adjoint eigenvalue problems, nonlinear eigenvalue problems and parameterized eigenvalue problems.

	
	\section*{Declarations}
	The authors have no competing interests to declare that are relevant to the content of this article.
	
	\bibliographystyle{siamplain}
	\bibliography{refs}
\end{document}